\documentclass[a4paper,12pt,reqno]{amsart}

\usepackage[headings]{fullpage}

\usepackage{amsfonts,graphics,amsmath,amsthm,amscd,amssymb,latexsym,euscript,enumerate,bm,booktabs}

\usepackage{epsfig}
\usepackage{flafter}
\usepackage[all,cmtip,line]{xy}
\usepackage{xltabular}
\usepackage{array}
\usepackage[english]{babel}
\usepackage{overpic}
\usepackage{subfig}
\usepackage{multirow}
\usepackage{microtype}
\usepackage{wrapfig}
\usepackage{tabularray}
\usepackage{longtable}
\usepackage{supertabular}
\usepackage{caption}

\usepackage[normalem]{ulem}
\usepackage{tikz}
\usetikzlibrary{positioning}
\usepackage{tikz-cd}
\usepackage{float}
\allowdisplaybreaks

\usepackage[dvipsnames,svgnames,table]{xcolor}
\usepackage{graphicx}
\usepackage{hyperref}
\hypersetup{
    colorlinks=true,    
    linkcolor=blue,          
    citecolor=blue,      
    filecolor=blue,      
    urlcolor=blue
}
\usepackage{makecell}

\newtheorem{theorem}{Theorem}[section]
\newtheorem*{maintheorem}{Main Theorem}
\makeatletter
\newcommand{\namedlabel}[2]{%
  \phantomsection
  \def\@currentlabel{#2}%
  \label{#1}%
}
\makeatother
\newtheorem{lemma}[theorem]{Lemma}
\newtheorem{proposition}[theorem]{Proposition}

\newtheorem{question}[theorem]{Question}
\newtheorem*{theorem*}{Theorem}

\theoremstyle{plain}
\newtheorem{corollary}[theorem]{Corollary}

\theoremstyle{definition} % italic or bold etc.

\newtheorem{definition-lemma}[theorem]{Definition-Lemma}

\newtheorem{remark}[theorem]{Remark}
\newtheorem*{notation}{Notation}
\numberwithin{equation}{section}

\newcommand{\C}{\mathbb{C}}
\newcommand{\R}{\mathbb{R}}
\newcommand{\Z}{\mathbb{Z}}

\newcommand{\Q}{\mathbb{Q}}
\newcommand{\OO}{\mathcal{O}}

\def\P{\mathbb{P}}

\newcommand{\A}{\mathbb{A}}
\newcommand{\F}{\mathbb{F}}

\DeclareMathOperator{\ord}{ord}

\DeclareMathOperator{\Pic}{Pic}

\def\mult{\operatorname{mult}}

\def\Supp{\operatorname{Supp}}
\def\Exc{\operatorname{Exc}}

\def\Sing{\operatorname{Sing}}

\usepackage{mathtools}

\DeclarePairedDelimiterX{\norm}[1]{\lVert}{\rVert}{#1}

\title[Cylindricity of log del Pezzo surfaces in characteristic five]
{Cylindricity of log del Pezzo surfaces in characteristic five}

\begin{document}

\author[D.-W. Lee]{Dae-Won Lee}
\author[M. Nagaoka]{Masaru Nagaoka}
\author[M. Sawahara]{Masatomo Sawahara}
\address[Dae-Won Lee]{School of Mathematics, Korea Institute for Advanced Study, 85 Hoegiro, Dongdaemun-gu, Seoul 02455, Republic of Korea}
\email{daewonlee@kias.re.kr}
\address[Masaru Nagaoka]{Gakushuin University, 1-5-1 Mejiro, Toshima-ku, Tokyo 171-8588,
Japan}
\email{masaru.nagaoka@gakushuin.ac.jp}
\address[Masatomo Sawahara]{Faculty of Education, Hirosaki University, Bunkyocho 1, Hirosaki-shi, Aomori 036-8560, Japan}
\email{sawahara.masatomo@gmail.com}

\subjclass[2020]{14J26, 14J45, 14G17, 14R25}
\date{\today}
\keywords{Log del Pezzo surface, cylinder, positive characteristic, log liftability}

\begin{abstract}
    Over an algebraically closed field of characteristic five, we prove that a log del Pezzo surface of Picard rank one contains no cylinder if it is not log liftable. We also characterize Du Val del Pezzo surfaces that contain no cylinder in positive characteristic.
\end{abstract}

\maketitle

\section{Introduction}\label{sect:1}

For a normal projective variety $X$, an open subset $U \subset X$ is called a \emph{cylinder} if $U$ is isomorphic to $Z \times \A^1$ for some affine variety $Z$. The existence of cylinders on Fano varieties is closely related to additive group actions on affine cones (see \cite{KPZ11}). In this paper, we study cylinders on log del Pezzo surfaces in positive characteristic.

Establishing the nonexistence of cylinders is often difficult. For surfaces of Picard rank one, Dubouloz and Kishimoto used Corti's inequality to establish the nonexistence of cylinders (see Lemma \ref{lem:2.3}). Indeed, Dubouloz and Kishimoto proved in \cite[Proposition 9]{DK18} that every smooth del Pezzo surface of Picard rank one and degree at most four defined over a field of characteristic zero contains no cylinder (see also {\cite[Theorem 3.30]{CPPZ}}). The result was extended to canonical del Pezzo surfaces in \cite{DKS}.

\begin{theorem}[{\cite{DKS}}]\label{thm:1.1}
Let $X$ be a canonical del Pezzo surface of Picard rank one defined over a field $\Bbbk$ with $\mathrm{char}(\Bbbk)=0$, and let $d$ be the degree of $X$; i.e., $d \coloneqq (-K_X)^2 \in \{1,\dots,6,8,9\}$. 
Assume that $X$ contains a cylinder. Then: 
\begin{enumerate}
\item If $d \le 4$, then $X$ has a singular point defined over $\Bbbk$. 
\item If $d=2$, then $X$ has a singular point defined over $\Bbbk$ not of type ${\sf A}_1$. 
\item If $d=1$, then $X$ has a singular point defined over $\Bbbk$ not of type ${\sf A}_1$, ${\sf A}_2$, ${\sf A}_3$, ${\sf D}_4$. 
\end{enumerate}
\end{theorem}

Combining Theorem \ref{thm:1.1} with the classification in \cite[Theorem 1.2 and Table 1.1]{Ye02} and the cylinder constructions in \cite[Section 4]{CPW16} and \cite[Theorem 1.6]{Bel17}, we obtain the following characterization: a Du Val del Pezzo surface $X$ of Picard rank one over an algebraically closed field of characteristic zero contains no cylinder if and only if its Dynkin type $\mathrm{Dyn}(X)$ is $4{\sf A}_2$, $2{\sf A}_1+2{\sf A}_3$, or $2{\sf D}_4$. The same characterization also follows by combining the classification of Du Val del Pezzo surfaces admitting anticanonically polar cylinders in \cite[Theorem 1.5]{CPW16} with \cite[Theorem 1.6]{Bel17}. For the cylindricity problem, however, the approach of Dubouloz and Kishimoto gives a shorter argument.

We first characterize Du Val del Pezzo surfaces that contain no cylinder in positive characteristic.

\begin{theorem}\label{thm:1.2}
Let $X$ be a Du Val del Pezzo surface defined over an algebraically closed field $\Bbbk$ with $\mathrm{char}(\Bbbk)>0$. Then $X$ contains no cylinder if and only if $\rho(X)=1$ and the Dynkin type ${\rm Dyn}(X)$ of $X$ is one of the following: 
\begin{align*}
4{\sf A}_2,\quad
2{\sf D}_4,\quad
2{\sf A}_1+2{\sf A}_3,\quad
4{\sf A}_1+{\sf D}_4,\quad
8{\sf A}_1,\quad
7{\sf A}_1. 
\end{align*}
\end{theorem}

We prove the nonexistence of cylinders using the local intersection inequalities in Lemma \ref{lem:2.3}, Corollary \ref{cor:2.6}, and Lemma \ref{lem:3.3}.

To prove the converse, Proposition \ref{prop:3.4} compares the configurations of $(-1)$-curves and $(-2)$-curves on minimal resolutions in different characteristics, and Proposition \ref{prop:3.5} establishes the existence of cylinders for the remaining Picard rank one surfaces. We then run the $K_X$-minimal model program. Lemmas \ref{lem:3.7} and \ref{lem:3.9} imply that the $K_X$-MMP of a Du Val del Pezzo surface $X$ containing no cylinder terminates with a Mori fiber space of Picard rank one, while Lemmas \ref{lem:3.8} and \ref{lem:3.6} exclude the possible nontrivial birational contractions.

We now relate the nonexistence of cylinders to log liftability over the ring of Witt vectors $W(\Bbbk)$. By \cite[Theorem 1.7]{KN22}, for every Du Val del Pezzo surface $X$ of Picard rank one defined over an algebraically closed field $\Bbbk$ with $\mathrm{char}(\Bbbk)>0$, $X$ is not log liftable over $W(\Bbbk)$ if and only if $(\mathrm{char}(\Bbbk),{\rm Dyn}(X)) = (3,4{\sf A}_2)$, $(2,4{\sf A}_1+{\sf D}_4)$, $(2,8{\sf A}_1)$ or $(2,7{\sf A}_1)$. Hence, Theorem \ref{thm:1.2} implies that every Du Val del Pezzo surface of Picard rank one which is not log liftable over $W(\Bbbk)$ contains no cylinder. Based on this observation, we propose the following question.

\begin{question}\label{ques:1.3}
Let $X$ be a log del Pezzo surface of Picard rank one defined over an algebraically closed field $\Bbbk$ with $\mathrm{char}(\Bbbk)>0$. If $X$ is not log liftable over $W(\Bbbk)$, does $X$ contain no cylinder?
\end{question}

Theorem \ref{thm:1.2} answers Question \ref{ques:1.3} affirmatively for Du Val del Pezzo surfaces. It remains to consider log del Pezzo surfaces with singularities worse than Du Val. The \ref{thm:main} gives a positive answer to Question \ref{ques:1.3} in characteristic five.

\begin{maintheorem}\namedlabel{thm:main}{Main Theorem}
Let $X$ be a log del Pezzo surface of Picard rank one over an algebraically closed field $\Bbbk$ of characteristic five. If $X$ is not log liftable over $W(\Bbbk)$, then $X$ contains no cylinder.
\end{maintheorem}

The proof of the \ref{thm:main} relies on the classification obtained by the second author \cite[Theorem 1.2 and Table 1]{Nag25}, summarized in Table \ref{tab:1}. A cylinder on a surface in this table determines a movable pencil with a unique base point at which the associated pair is not log canonical. Corollary \ref{cor:2.6} excludes smooth base points and the singular points whose exceptional divisors have sufficiently small determinants. This proves the assertion for rows I--III and Nos. 1--3.

For Nos. 4--21, these local bounds alone do not exclude every singular base point. The curve constructed in Subsection \ref{subsec:4.3} meets exactly one component of the exceptional divisor over each remaining candidate base point, transversely at a single point. This allows us to bound the coefficients along the exceptional curves in the pullback of a general member of the pencil. Proposition \ref{prop:2.8} combines these bounds with Lemma \ref{lem:2.2} and Corti's inequality. The determinant and discrepancy computations in Section \ref{sect:4} verify the criterion for each remaining candidate base point.

Our argument is closely related to the method of maximal singularities in birational rigidity (see \cite[Section 3]{Cor00}). In our setting, Lemma \ref{lem:2.1} gives a divisor with discrepancy $-2$, and this failure of log canonicity gives a lower bound on local intersection multiplicities that is incompatible with the corresponding global intersection number. Proposition \ref{prop:2.8} formulates this argument as a criterion on the minimal resolution.

The rest of this paper is organized as follows. Section \ref{sect:2} contains the results on pencils associated with a cylinder and the local intersection inequalities. Section \ref{sect:3} proves Theorem \ref{thm:1.2}, and Section \ref{sect:4} proves the \ref{thm:main} using the constructions in Subsection \ref{subsec:4.3} and the exclusion arguments in Propositions \ref{prop:4.5} and \ref{prop:4.9}.

\section{Preliminaries}\label{sect:2}
\subsection{Notation and conventions}\label{subsec:2.1}
Unless otherwise stated, all varieties are defined over an algebraically closed field $\Bbbk$. The results of Sections \ref{sect:2} and \ref{sect:3} hold in arbitrary characteristic, whereas throughout Section \ref{sect:4} we assume $\mathrm{char}(\Bbbk)=5$. A \emph{log del Pezzo surface} is a normal projective surface with klt singularities and ample anticanonical divisor. A \emph{Du Val del Pezzo surface} is a log del Pezzo surface with only Du Val singularities. When $\mathrm{char}(\Bbbk)>0$, a normal projective surface $S$ is \emph{log liftable over $W(\Bbbk)$} if, for some log resolution $f\colon V\to S$, the pair $(V,\Exc(f)_{\rm red})$ is the special fiber of a smooth projective scheme over $W(\Bbbk)$ equipped with a relative simple normal crossings divisor (see \cite[Definitions 3.1--3.2]{Nag25}).

\begin{notation}
    Throughout this paper, the rank of a surface means its Picard rank. Let $A=A_1+\cdots+A_s$ be an ordered chain of smooth rational curves with simple normal crossings. We say that $A$ has \emph{type} $[a_1,\dots,a_s]$ when $A_i^2=-a_i$, $A_i\cdot A_{i+1}=1$, and $A_i\cdot A_j=0$ for $|i-j|>1$. Unless explicitly stated otherwise, the integers $a_i$ are at least two. The type $[b;B_1,B_2,B_3]$ describes a reduced divisor consisting of a central curve of self-intersection $-b$ and three chains with the given types. Each chain is ordered starting from the component meeting the central curve. The three chains are pairwise disjoint, and the central curve meets only the first component of each chain, transversely once. The notation $2^t$ denotes $t$ consecutive entries equal to two.

    For an ordered reduced divisor $E=E_1+\cdots+E_s$ on a smooth surface, let $C_E\coloneqq (-E_i\cdot E_j)_{i,j}$ and $d(E)\coloneqq \det C_E$. If $C_E$ is invertible, then set $G_E\coloneqq C_E^{-1}$. For a chain $A$ of type $[a_1,\dots,a_s]$, we also write $d([a_1,\dots,a_s])$ for $d(A)$. We set $d(0)=d([])=1$. When a reduced exceptional divisor $E$ is fixed, we simply write $C_E,G_E,d(E)$ as $C,G,D$, respectively. For index sets $I,J$, the notation $G_{I,J}$ denotes the submatrix with rows in $I$ and columns in $J$. 

    For a normal $\Q$-Gorenstein surface $S$, let $\pi\colon Y\to S$ be the minimal resolution, and set $H\coloneqq \pi^{\ast}(-K_S)$. When $S$ is projective, set $h\coloneqq H^2=(-K_S)^2$.
    For every exceptional prime divisor $T\subset Y$, the coefficient $\delta_T$ and the number $\alpha_T$ are defined by
    \[
    \pi^{\ast}K_S=K_Y+\sum_{T\subset \Exc(\pi)}\delta_T T,\quad \alpha_T\coloneqq 1-\delta_T.
    \]
    If $\Exc(\pi)_{\rm red}=E_1+\cdots+E_s$, we write $\delta_i\coloneqq\delta_{E_i}$ and $\alpha_i\coloneqq\alpha_{E_i}$ for $1\leq i\leq s$. A tilde denotes the strict transform on $Y$. When a curve $L_Y\subset Y$ not contracted by $\pi$ is specified, we let $L\coloneqq \pi_{\ast}L_Y$ and $\eta\coloneqq H\cdot L_Y=(-K_S)\cdot L$. 
\end{notation}

\subsection{Pencils associated to cylinders}\label{subsec:2.2}
Let $\mathcal{M}$ be a movable linear system on a normal surface $S$, let $B$ be a $\Q$-divisor, and let $c\in \Q_{\geq 0}$. Assume that $K_S+B$ and every member of $\mathcal{M}$ are $\Q$-Cartier. For a prime divisor $F$ over $S$, choose a resolution $\mu\colon W\to S$ on which $F$ appears as a prime divisor, and choose a general member $M\in \mathcal{M}$. For an exceptional divisor $F$, the \emph{discrepancy} $a_F(S,c\mathcal{M}+B)$ is the coefficient of $F$ in
\[
K_W+cM_W+B_W-\mu^{\ast}(K_S+cM+B),
\]
where $M_W$ and $B_W$ are the strict transforms of $M$ and $B$, respectively. For a prime divisor $F\subset S$, the discrepancy is defined as $a_F(S,c\mathcal{M}+B)\coloneqq -\mult_F(B)$, since $\mathcal{M}$ has no fixed component. We say that $(S,c\mathcal{M}+B)$ is \emph{log canonical} if $a_F(S,c\mathcal{M}+B)\geq-1$ for every prime divisor $F$ over $S$. 

For a morphism $f\colon S\to C$ from a surface to a curve, an irreducible curve $D\subset S$ is \emph{vertical} if $f(D)$ is a point and \emph{horizontal} if $f(D)=C$.

\begin{lemma}\label{lem:2.1}
    Let $S$ be a normal projective $\Q$-factorial rational surface of Picard rank one with $-K_S$ ample. Suppose that $S$ contains a cylinder $U\simeq \A^1\times Z$. Then the closures of the fibers of the projection $U\to Z$ form a movable pencil $\mathcal{L}$ on $S$, and there is a rational number $a>0$ such that $\mathcal{L}\equiv -aK_S$. Moreover, the pencil $\mathcal{L}$ has a unique base point $P$, and the pair $(S,\frac{1}{a}\mathcal{L})$ is not log canonical at $P$.

    For the strict transform $\widetilde{\mathcal{L}}$ on the minimal resolution $\pi\colon Y\to S$, the following assertions hold:
    \begin{enumerate}[(1)]
        \item If $\widetilde{\mathcal{L}}^2=0$, then $\widetilde{\mathcal{L}}$ has no base points, and exactly one irreducible component of $\pi^{-1}(P)$ is horizontal for the morphism $Y\to \P^1$ defined by $\widetilde{\mathcal{L}}$, and that component is a section.
        \item If $\widetilde{\mathcal{L}}^2>0$, then $\widetilde{\mathcal{L}}$ has a unique base point $q\in \pi^{-1}(P)$. Let $\psi\colon W\to Y$ be a minimal sequence of point blow-ups resolving the base points of $\widetilde{\mathcal{L}}$, and let $f\colon W\to \P^1$ be the induced morphism. Then the exceptional curve of the last blow-up is a section of $f$, and all other curves contracted by $\psi$, together with the strict transforms of the curves contracted by $\pi$, are vertical for $f$.
    \end{enumerate}
\end{lemma}
\begin{proof}
    Since $\A^1\times Z$ is normal, the curve $Z$ is normal. Since $\Bbbk$ is perfect, $Z$ and $U$ are smooth. In particular, $U$ is contained in the smooth locus of $S$. Let $\overline{Z}$ be the smooth projective completion of $Z$. Since $S$ is rational, we obtain a dominant rational map $\P^2\dashrightarrow \overline{Z}$. By restricting the rational map to a general line, we obtain a nonconstant rational map $\P^1\dashrightarrow \overline{Z}$, and L\"uroth's theorem implies that $\overline{Z}\simeq \P^1$. The projection $U\to Z$ extends to a rational map $S\dashrightarrow \overline{Z}\simeq \P^1$, which defines a pencil $\mathcal{L}$ on $S$, and a general member of $\mathcal{L}$ is integral and meets $U$ in an affine line. Moreover, $\mathcal{L}$ has no fixed components.

    Since $\rho(S)=1$ and $-K_S$ is ample, a general member of $\mathcal{L}$ is numerically equivalent to $-aK_S$ for some $a\in \Q_{>0}$. Hence, two distinct general members $L_1$ and $L_2$ of $\mathcal{L}$ satisfy $L_1\cdot L_2=a^2(-K_S)^2>0$. Since fibers of a morphism to a curve are disjoint, the pencil $\mathcal{L}$ must have a base point.

    Choose a resolution $\mu\colon W\to S$ of the pencil which is an isomorphism over $U$, with $W$ smooth and $W\setminus U$ a simple normal crossings divisor. Write $f\colon W\to \P^1$ for the induced morphism. Then the generic fiber of $f$ is the regular projective completion of $\A^1_{\Bbbk(Z)}$ and hence is $\P^1_{\Bbbk(Z)}$. The complement of $U$ in the generic fiber consists of the single point at infinity defined over $\Bbbk(Z)$. Therefore, $W\setminus U$ has exactly one horizontal component $F_{\rm hor}$, and the degree of the map $F_{\rm hor}\to \P^1$ is one.

    Write $\Phi\colon S\dashrightarrow \P^1$ for the map defined by $\mathcal{L}$, and let $V\subset S\times \P^1$ be the closure of $\{(x,\Phi(x))\mid x\in \mathrm{Dom}(\Phi)\}$. Denote the projections by $\nu\colon V\to S$ and $g\colon V\to \P^1$. If the fiber of $\nu$ over a point has dimension zero, then properness implies that $\nu$ is finite in a neighborhood of that point. Since $S$ is normal, every finite birational morphism to $S$ is an isomorphism. Thus, a base point of $\mathcal{L}$ has a fiber of dimension one under $\nu$. Such a fiber contains a curve dominating $\P^1$ under $g$. Hence, a resolution of $V$ has a horizontal exceptional curve above each base point. Uniqueness of $F_{\rm hor}$ now gives a unique base point $P$, with $\mu(F_{\rm hor})=P$. The smooth rational curve $F_{\rm hor}$ maps with degree one to $\P^1$, hence is a section of $f$.

    Let $F$ be a general fiber of $f$, and let $L_F\coloneqq \mu_{\ast}F\in \mathcal{L}$. Denote the curves contracted by $\mu$ by $F_1,\dots,F_t=F_{\rm hor}$. Since $F$ is the strict transform of $L_F$, we have
    \begin{equation}\label{eq:2.1}
        K_W+\frac{1}{a}F=\mu^{\ast}\left(K_S+\frac{1}{a}L_F\right)+\sum_{i=1}^t e_iF_i.
    \end{equation}
    For $i<t$, we have $F\cdot F_i=0$ and $F\cdot F_t=1$. Since $K_S+\frac{1}{a}L_F\equiv 0$ and $K_W\cdot F=-2$, \eqref{eq:2.1} gives $e_t=-2$. Thus, $F_t$ has discrepancy less than $-1$, which implies that $(S,\frac{1}{a}\mathcal{L})$ is not log canonical at $P$.

    On the minimal resolution $Y$, the pencil $\widetilde{\mathcal{L}}$ has nonnegative self-intersection number. Resolve the base points of $\widetilde{\mathcal{L}}$ by a sequence of $v$ point blow-ups, and let $m_i$ be the multiplicity of the strict transform of a general member at the center of the $i$th blow-up, for $1\leq i\leq v$. The self-intersection number drops by $m_i^2$ at the $i$th blow-up, and the final pencil has self-intersection zero. Hence, $\widetilde{\mathcal{L}}^2=\sum_{i=1}^v m_i^2$, which vanishes exactly when $\widetilde{\mathcal{L}}$ has no base points. If $\widetilde{\mathcal{L}}$ has no base points, then some component of $\pi^{-1}(P)$ is horizontal. By the uniqueness of the horizontal boundary component proved above, it is the only such component and has degree one over $\P^1$. Hence, it is a section, proving (1).

    If $\widetilde{\mathcal{L}}^2>0$, then every base point of $\widetilde{\mathcal{L}}$ has a horizontal exceptional curve above the base point. Thus, $\widetilde{\mathcal{L}}$ has a unique base point $q$ on $Y$ lying over $P$. In the resolution $\psi$, the intersection number of the last exceptional curve with a general fiber equals the multiplicity of the pencil at the center of the last blow-up. Therefore, the last exceptional curve is the only horizontal component of $W\setminus U$, and its degree over $\P^1$ is one. All other exceptional curves on $W$ are vertical, proving (2).
\end{proof}

\begin{lemma}\label{lem:2.2}
    Let $\mathcal{L}$ be a pencil as in Lemma \ref{lem:2.1}, with base point $P\in S$. Suppose that $\widetilde{\mathcal{L}}^2>0$, and denote the base point of $\widetilde{\mathcal{L}}$ by $q$. Let $\Theta\geq 0$ be a $\Q$-divisor supported on $\pi^{-1}(P)$, and let $c\in \Q_{>0}$. If
    \begin{equation}\label{eq:2.2}
        \widetilde{\mathcal{L}}\cdot (K_Y+c\widetilde{\mathcal{L}}+\Theta)\geq 0,
    \end{equation}
    then the pair $(Y,c\widetilde{\mathcal{L}}+\Theta)$ is not log canonical at $q$.
\end{lemma}
\begin{proof}
    Let $\psi\colon W\to Y$ be a minimal sequence of point blow-ups resolving the base points of $\widetilde{\mathcal{L}}$, let $f\colon W\to \P^1$ be the induced morphism, and let $F$ be a general fiber of $f$. Let $M_F\coloneqq \psi_{\ast}F\in \widetilde{\mathcal{L}}$, and denote the strict transform of $\Theta$ by $\Theta_W$. Label the curves contracted by $\psi$ as $F_1,\dots,F_t$ in the order of the blow-ups. Then we have
    \[
    K_W+cF+\Theta_W=\psi^{\ast}(K_Y+cM_F+\Theta)+\sum_{i=1}^t e_iF_i.
    \]
    By Lemma \ref{lem:2.1}(2), every component of $\Theta_W$ and each $F_i$ for $i<t$ are vertical for $f$, whereas $F_t$ is a section. Since we have $-2=\widetilde{\mathcal{L}}\cdot (K_Y+c\widetilde{\mathcal{L}}+\Theta)+e_t$, \eqref{eq:2.2} implies that $e_t\leq -2<-1$. Thus, $(Y,c\widetilde{\mathcal{L}}+\Theta)$ is not log canonical at $q$.
\end{proof}

\subsection{Corti's inequality on a smooth surface}
Corti's inequality was established in characteristic zero in \cite[Theorem 3.1]{Cor00}. The proof below applies in arbitrary characteristic.
\begin{lemma}[cf. \protect{\cite[Theorem 3.1]{Cor00}}]\label{lem:2.3}
    Let $p$ be a smooth point of a surface $S$, and let $B_1,B_2\subset S$ be smooth curves meeting transversely at $p$. Let $\mathcal{M}$ be a movable linear system, let $c\in \Q_{>0}$, and let $b_1,b_2\in[0,1]\cap \Q$. If $(S,c\mathcal{M}+b_1B_1+b_2B_2)$ is not log canonical at $p$, then distinct general members $M_1,M_2\in \mathcal{M}$ satisfy
    \[
    c^2(M_1\cdot M_2)_{p}>4(1-b_1)(1-b_2).
    \]
\end{lemma}
\begin{proof}
Put $\varepsilon_j\coloneqq1-b_j$ for $j=1,2$, $m\coloneqq c\mult_p\mathcal{M}$, and $I_p\coloneqq c^2(M_1\cdot M_2)_p$. Choose an exceptional prime divisor $F$ with center $p$ and discrepancy less than $-1$. Such a divisor exists since $\mathcal{M}$ has no fixed component and $b_1,b_2\leq1$. We proceed by induction on the number of point blow-ups required to extract $F$.

Let $\sigma\colon S'\to S$ be the blow-up of $p$, with exceptional curve $E$, and write $\mathcal{M}',B_1',B_2'$ for the strict transforms. Let $M'$ be the strict transform of a general member $M\in\mathcal M$. The equality
\[
 K_{S'}+cM'+b_1B_1'+b_2B_2'+(1-\varepsilon_E)E
 =\sigma^*(K_S+cM+b_1B_1+b_2B_2)
\]
gives $\varepsilon_E\coloneqq\varepsilon_1+\varepsilon_2-m$. If $\varepsilon_E<0$, then $I_p\geq m^2>(\varepsilon_1+\varepsilon_2)^2\geq4\varepsilon_1\varepsilon_2$. The assertion follows whenever $\varepsilon_E<0$, including the case $F=E$. We may therefore assume $\varepsilon_E\geq0$, and the center of $F$ on $S'$ is a point $q\in E$. For the strict transforms $M_1',M_2'$, put $I_q\coloneqq c^2(M_1'\cdot M_2')_q$. The intersection formula for $\sigma$ gives $I_p\geq m^2+I_q$. The curves $B_1'$ and $B_2'$ meet $E$ at distinct points.

If $\varepsilon_E>1$, then replace the negative coefficient $1-\varepsilon_E$ of $E$ by zero. The discrepancy of $F$ does not increase. Induction gives $I_q>4\varepsilon_j$ when $q\in B_j'$ and $I_q>4$ when $q\notin B_1'\cup B_2'$. Each inequality implies $I_p>4\varepsilon_1\varepsilon_2$, since $0\leq\varepsilon_1,\varepsilon_2\leq1$.

Suppose that $0\leq\varepsilon_E\leq1$. If $q\in B_j'\cap E$ for $j\in\{1,2\}$, then induction gives $I_q>4\varepsilon_j\varepsilon_E$. Consequently, $I_p>m^2+4\varepsilon_j(\varepsilon_1+\varepsilon_2-m)=(m-2\varepsilon_j)^2+4\varepsilon_1\varepsilon_2\geq4\varepsilon_1\varepsilon_2$. If $q\in E\setminus(B_1'\cup B_2')$, then induction gives $I_q>4\varepsilon_E$. Put $t\coloneqq\varepsilon_1+\varepsilon_2$. Since $0\leq m\leq t\leq2$, we have $m^2+4(t-m)-t^2=(t-m)(4-t-m)\geq0$. Hence, $I_p>m^2+4(t-m)\geq t^2\geq4\varepsilon_1\varepsilon_2$. If only one boundary component passes through the center of $F$, choose another smooth curve through that point which meets the first one transversely, and assign coefficient zero to it. Thus, the induction hypothesis can be applied in each case.
\end{proof}

\subsection{Intersection matrices of chains}
Let $A=A_1+\cdots+A_s$ be an ordered chain of type $[a_1,\dots,a_s]$, and set $C\coloneqq C_A$ and $D\coloneqq d(A)$. For $1\leq i\leq s$, define
\begin{equation}\label{eq:2.3}
    p_i\coloneqq d([a_1,\dots,a_{i-1}]),\quad q_i\coloneqq d([a_{i+1},\dots,a_s]).
\end{equation}
The convention $d(0)=1$ gives $p_1=q_s=1$. We also set $p_0=q_{s+1}=0$ and $p_{s+1}=q_0=D$.

\begin{lemma}\label{lem:2.4}
    For $1\leq i\leq s$, the determinants in \eqref{eq:2.3} satisfy
    \begin{equation}\label{eq:2.4}
        p_{i+1}=a_ip_i-p_{i-1},\quad q_{i-1}=a_iq_i-q_{i+1}.
    \end{equation}
    The sequence $p_0,\dots,p_{s+1}$ is strictly increasing, and $q_0,\dots,q_{s+1}$ is strictly decreasing. Moreover, we have
    \begin{equation}\label{eq:2.5}
    \begin{aligned}
    (C^{-1})_{ij}&=\frac{p_{\min\{i,j\}}q_{\max\{i,j\}}}{D}
    &&(1\leq i,j\leq s),\\
    p_{i+1}q_i-p_iq_{i+1}&=D &&(1\leq i\leq s).
    \end{aligned}
    \end{equation}
    If $A$ is the reduced exceptional divisor over a point of a normal $\Q$-Gorenstein surface $S$, then 
    \begin{equation}\label{eq:2.6}
        \delta_i=1-\frac{p_i+q_i}{D},\quad 0\leq\delta_i<1.
    \end{equation}
\end{lemma}
\begin{proof}
The recurrences \eqref{eq:2.4} agree with those in \cite[Section 1, equation (2)]{Sug12}, with diagonal entries $a_i$ and entries $-1$ on both adjacent diagonals. Since $a_i\geq2$, induction gives $p_{i+1}-p_i>0$ and $q_i-q_{i+1}>0$ from the initial differences $p_1-p_0=q_s-q_{s+1}=1$. The leading principal minors $p_2,\ldots,p_{s+1}$ are therefore positive. Hence, $C$ is positive definite.

The formula for $C^{-1}$ in \eqref{eq:2.5} is given in \cite[Section 1, equation (2)]{Sug12}. The determinants denoted there by $\theta_{i-1}$, $\phi_{i+1}$, and $\theta_s$ correspond to $p_i$, $q_i$, and $D$, respectively. The recurrences \eqref{eq:2.4} show that $p_{i+1}q_i-p_iq_{i+1}$ is independent of $i$; at $i=s$ its value is $D$.

Let $\bm{\delta}=(\delta_1,\dots,\delta_s)^{\rm t}$ and $\bm{\alpha}=(\alpha_1,\dots,\alpha_s)^{\rm t}$. By adjunction, we obtain $C\bm{\delta}=(a_1-2,\dots,a_s-2)^{\rm t}$ and $C\bm{\alpha}=\bm{e}_1+\bm{e}_s$, where $\bm{e}_i$ is the $i$th coordinate vector. By \eqref{eq:2.5}, we obtain $\alpha_i=(C^{-1})_{i1}+(C^{-1})_{is}=\frac{p_i+q_i}{D}$. Since $C^{-1}$ has positive entries and $a_i-2\geq 0$, we have $\delta_i\geq 0$. The inequality $\alpha_i>0$ gives $\delta_i<1$, which proves \eqref{eq:2.6}.
\end{proof}

\begin{lemma}\label{lem:2.5}
    Let $P$ be a singular point of a normal $\Q$-Gorenstein surface $S$ whose reduced exceptional divisor on the minimal resolution is a chain $A=A_1+\cdots+A_s$ of type $[a_1,\dots,a_s]$, with $a_i\geq 2$. Let $\mathcal{M}$ be a movable linear system whose members are $\Q$-Cartier, and let $c\in \Q_{>0}$. If $(S,c\mathcal{M})$ is not log canonical at $P$, then for distinct general members $M_1,M_2\in \mathcal{M}$, we have
    \begin{equation}\label{eq:2.7}
        c^2(M_1\cdot M_2)_P>\frac{4}{d(A)}.
    \end{equation}
\end{lemma}
\begin{proof}
    In the following computation, we consider only the exceptional components lying over $P$. For a general member $M\in \mathcal{M}$, write $c\pi^{\ast}M=c\widetilde{M}+\sum_{i=1}^s x_iA_i$ and $\bm{x}=(x_1,\dots,x_s)^{\rm t}$. Since $\pi^{\ast}M\cdot A_i=0$ for every $i$, we have $C\bm{x}=(c\widetilde{M}\cdot A_1,\dots, c\widetilde{M}\cdot A_s)^{\rm t}$. Note that the entries of $C\bm{x}$ are nonnegative, and hence $x_i\geq 0$. By \eqref{eq:2.6}, we have $\alpha_i=\frac{p_i+q_i}{D}$.

    Choose two distinct general members $M_1,M_2\in\mathcal M$ whose pullbacks have the same coefficients along each $A_i$, and let $\widetilde M_1$ and $\widetilde M_2$ be their strict transforms, respectively. Then we have
    \begin{equation}\label{eq:2.8}
        c^2(M_1\cdot M_2)_P=\bm{x}^{\rm t}C\bm{x}+\sum_{q\in \pi^{-1}(P)}c^2(\widetilde{M}_1\cdot \widetilde{M}_2)_q.
    \end{equation}

    For a general member $M\in\mathcal M$, we have, in a neighborhood of $\pi^{-1}(P)$,
    \[
    \pi^*(K_S+cM)=K_Y+c\widetilde M+\sum_{i=1}^s(\delta_i+x_i)A_i.
    \] 
    If $x_i>\alpha_i$ for some $i$, then by the Cauchy--Schwarz inequality, we have
    \[
    \bm{x}^{\rm t}C\bm{x}\geq \frac{x_i^2}{G_{ii}}>\frac{\alpha_i^2}{G_{ii}}=\frac{(p_i+q_i)^2}{Dp_iq_i}\geq \frac{4}{D}.
    \]
    This proves \eqref{eq:2.7} in this case. We may now assume $0\leq x_i\leq \alpha_i$ for every $i$. Since $(S,c\mathcal M)$ is not log canonical at $P$, the pair
    \[
    \left(Y,c\widetilde{\mathcal M}+\sum_{i=1}^s(\delta_i+x_i)A_i\right)
    \]
    is not log canonical at some point of $A$.

    Suppose first that the pair is not log canonical at a point of $A_i\setminus\bigcup_{j\neq i}A_j$. The coefficient of $A_i$ is $1-(\alpha_i-x_i)\in [0,1]$. By Lemma \ref{lem:2.3} and \eqref{eq:2.8}, we have
    \[
    c^2(M_1\cdot M_2)_P>\frac{x_i^2}{G_{ii}}+4(\alpha_i-x_i).
    \]
    For $0\leq z\leq \alpha_i$, the derivative of $\frac{z^2}{G_{ii}}+4(\alpha_i-z)$ is nonpositive, since $\alpha_i\leq 2G_{ii}$ is equivalent to $p_i+q_i\leq 2p_iq_i$. Thus, $\frac{z^2}{G_{ii}}+4(\alpha_i-z)\geq \frac{\alpha_i^2}{G_{ii}}\geq \frac{4}{D}$.

    Suppose now that the pair is not log canonical at $A_i\cap A_{i+1}$. Let $I\coloneqq \{i,i+1\}$, $z\coloneqq x_i$, and $w\coloneqq x_{i+1}$. By \eqref{eq:2.5},
    \[
    (G_{I,I})^{-1}=
    \begin{pmatrix}
    \frac{p_{i+1}}{p_i}&-1\\-1&\frac{q_i}{q_{i+1}}
    \end{pmatrix}.
    \]
    Let $\bm{v}\coloneqq (z,w)^{\rm t}$ and $\bm{y}\coloneqq G_{\{1,\dots,s\},I}(G_{I,I})^{-1}\bm{v}$. Then $y_i=z$, $y_{i+1}=w$, and every entry of $C\bm{y}$ outside $I$ is zero. Hence, $(\bm{x}-\bm{y})^{\rm t}C\bm{y}=0$, and 
    \[
    \bm{x}^{\rm t}C\bm{x} \geq \frac{p_{i+1}}{p_i}z^2 -2zw + \frac{q_i}{q_{i+1}}w^2.
    \]
    By Lemma \ref{lem:2.3}, $c^2(\widetilde{M}_1\cdot\widetilde{M}_2)_{A_i\cap A_{i+1}} > 4(\alpha_i-z)(\alpha_{i+1}-w)$. Thus, the sum of the two lower bounds is
    \begin{align*}
    &\frac{p_{i+1}}{p_i} z^2-2zw+\frac{q_i}{q_{i+1}}w^2
    +4\left(\frac{p_i+q_i}{D}-z\right)\left(\frac{p_{i+1}+q_{i+1}}D-w\right)\\
    &\quad=\frac4D+
    \begin{pmatrix}
    z-\frac{2p_i}{D}&w-\frac{2q_{i+1}}{D}
    \end{pmatrix}
    \begin{pmatrix}
    \frac{p_{i+1}}{p_i}&1\\1&\frac{q_i}{q_{i+1}}
    \end{pmatrix}
    \begin{pmatrix}
    z-\frac{2p_i}{D}\\w-\frac{2q_{i+1}}{D}
    \end{pmatrix}.
    \end{align*}
    The matrix
    \[
    \begin{pmatrix}
    \frac{p_{i+1}}{p_i}&1\\1&\frac{q_i}{q_{i+1}}
    \end{pmatrix}
    \]
    is positive definite, since the first diagonal entry is positive and the determinant is $\frac{D}{p_iq_{i+1}}>0$. Hence, the value is at least $\frac{4}{D}$. This completes the proof.
    \end{proof}

\begin{corollary}\label{cor:2.6}
    Let $\mathcal{L}\equiv -aK_S$ be the pencil associated to a cylinder as in Lemma \ref{lem:2.1}, and let $P$ be the base point. If $P$ is smooth, then $(-K_S)^2>4$. If the reduced exceptional divisor on the minimal resolution over $P$ is a chain $A$, then $(-K_S)^2>\frac{4}{d(A)}$.
\end{corollary}
\begin{proof}
    For distinct general members $L_1,L_2\in \mathcal{L}$, we have $L_1\cdot L_2=a^2(-K_S)^2$, and $(L_1\cdot L_2)_P\leq L_1\cdot L_2$. By Lemma \ref{lem:2.3} at a smooth point $P$, or Lemma \ref{lem:2.5} at a point with chain $A$, the assertion follows.
\end{proof}

\begin{remark}
When $S$ is smooth and $\mathrm{char}(\Bbbk)=0$, the conclusion of Corollary \ref{cor:2.6} also follows from \cite[Proposition 9]{DK18}. Our proof applies in arbitrary characteristic and also treats singular base points.
\end{remark}

\subsection{A criterion for excluding a base point}
\begin{lemma}\label{lem:2.7}
    Let $C$ be a real symmetric positive definite $s\times s$ matrix. Assume that $C_{ij}\leq 0$ for $i\neq j$ and that no simultaneous permutation of rows and columns puts $C$ into block diagonal form with two or more nonempty blocks. Let $G\coloneqq C^{-1}$. Then every entry of $G$ is positive, and
    \begin{equation}\label{eq:2.9}
        \frac{G_{ij}}{G_{Rj}}\leq \frac{G_{ii}}{G_{Ri}}\quad (1\leq R,i,j\leq s).
    \end{equation}
\end{lemma}
\begin{proof}
    Choose a real number $\lambda$ at least the largest eigenvalue of $C$, and put $M_0\coloneqq\lambda\mathrm{I}_s-C$. Every entry of $M_0$ is nonnegative, and every eigenvalue of $M_0$ lies in $[0,\lambda)$. Thus,
\[
 G=\frac1\lambda\sum_{n=0}^{\infty}\left(\frac{M_0}{\lambda}\right)^n
\]
has nonnegative entries. The assumption on block decompositions implies that, for distinct $i,j$, there are indices $i=i_0,i_1,\ldots,i_t=j$ with $(M_0)_{i_{v-1},i_v}>0$ for $1\leq v\leq t$. Hence, $(M_0^t)_{ij}>0$ and $G_{ij}>0$. The term $n=0$ gives $G_{ii}>0$. If $R=i$ or $j=i$, equality holds in \eqref{eq:2.9}. Otherwise, set $J=\{1,\dots,s\}\setminus \{i\}$. Then the same argument shows that $(C_{J,J})^{-1}$ has nonnegative entries. Since
\[
(C_{J,J})^{-1}=G_{J,J}-\frac{G_{J,\{i\}}G_{\{i\},J}}{G_{ii}},
\]
we obtain $G_{Rj}\geq \frac{G_{Ri}G_{ij}}{G_{ii}}$, which proves \eqref{eq:2.9}.
\end{proof}

\begin{proposition}\label{prop:2.8}
    Let $S$ be a $\Q$-factorial rational log del Pezzo surface of Picard rank one. Let $\pi\colon Y\to S$ be the minimal resolution. Let $p\in S$ be a singular point, and write the reduced exceptional divisor over $p$ as $E\coloneqq\pi^{-1}(p)_{\rm red}=E_1+\cdots+E_s$, where the $E_i$ are its irreducible components. Assume that $E$ has simple normal crossings. Fix $R\in \{1,\dots,s\}$, and suppose that an irreducible curve $L_Y\subset Y$, not contracted by $\pi$, meets $E_R$ transversely at a single point and is disjoint from $E_i$ for $i\neq R$.

    Let $C,G,h,\eta$ and $\delta_i$ be as in Subsection \ref{subsec:2.1} for $E$ and $L_Y$, and assume $\delta_i\geq 0$. For the curve $L_Y$, define
    \begin{equation*}
    \begin{aligned}
    \gamma_p\coloneqq 1+\frac{\eta}{h}\max_{1\leq i\leq s}\frac{\delta_i}{G_{Ri}},\quad
    \beta_p\coloneqq \gamma_p\eta \max_{1\leq i\leq s}\frac{G_{ii}}{G_{Ri}}.
    \end{aligned}
    \end{equation*}
    If 
    \begin{equation}\label{eq:2.10}
        \beta_p<1\quad \text{and}\quad \gamma_p^2h\leq 4(1-\beta_p)^2,
    \end{equation}
    then $p$ is not the base point of any pencil associated to a cylinder in $S$.
\end{proposition}
\begin{proof}
    Suppose that a cylinder in $S$ has associated pencil $\mathcal{L}\equiv-aK_S$ with base point $p$. A general member of $\mathcal{L}$ does not pass through any singular point of $S$ other than $p$. Therefore,
\[
 \widetilde{\mathcal{L}}\equiv aH-a\sum_{i=1}^s z_iE_i
\]
for rational numbers $z_i$. Put $r_i\coloneqq\frac{\widetilde{\mathcal{L}}\cdot E_i}{a}$, $\bm{r}\coloneqq(r_1,\ldots,r_s)^{\rm t}$, and $\bm{z}\coloneqq(z_1,\ldots,z_s)^{\rm t}$. Since $\widetilde{\mathcal{L}}$ has no fixed component, $r_i\geq0$ and $\widetilde{\mathcal{L}}\cdot L_Y\geq0$. Moreover, since
\[
\widetilde{\mathcal{L}}\cdot E_i=ar_i \quad\text{and}\quad \widetilde{\mathcal{L}}\cdot L_Y=a(\eta-z_R)\geq0,
\]
we have
\begin{equation}\label{eq:2.11}
\bm{z}=G\bm{r},\qquad z_R\leq\eta,\qquad
\widetilde{\mathcal{L}}^2=a^2(h-\bm{r}^{\rm t}G\bm{r}).
\end{equation} 
The divisor $E$ is connected and has negative definite intersection matrix. Thus, $C$ satisfies the assumptions of Lemma \ref{lem:2.7}, and $z_i\geq0$ for every $i$.

Define $\Theta\coloneqq\gamma_p\sum_{i=1}^s z_iE_i$. By Lemma \ref{lem:2.7}, the coefficient of $E_i$ in $\Theta$ satisfies
\[
 \gamma_pz_i
 \leq\gamma_p\frac{G_{ii}}{G_{Ri}}\sum_{j=1}^sG_{Rj}r_j
 \leq\gamma_p\eta\frac{G_{ii}}{G_{Ri}}
 \leq\beta_p.
\]
The definition of $\gamma_p$ also gives
\[
 \sum_{i=1}^s\delta_ir_i
 \leq(\gamma_p-1)\frac h\eta\sum_{i=1}^sG_{Ri}r_i
 \leq(\gamma_p-1)h.
\]
A general member of $\widetilde{\mathcal{L}}$ is disjoint from every exceptional curve over $\Sing(S)\setminus\{p\}$. From $K_Y=\pi^*K_S-\sum_{T\subset \Exc(\pi)}\delta_T T$, we obtain $\widetilde{\mathcal L}\cdot K_Y=-a(h+\sum_{i=1}^s\delta_i r_i)$, while $\widetilde{\mathcal{L}}\cdot\Theta=\gamma_pa\bm{r}^{\rm t}G\bm{r}$. Together with \eqref{eq:2.11}, these equalities yield
\begin{equation}\label{eq:2.12}
 \widetilde{\mathcal{L}}\cdot\left(K_Y+\frac{\gamma_p}{a}\widetilde{\mathcal{L}}+\Theta\right)
 =a\left((\gamma_p-1)h-\sum_{i=1}^s\delta_ir_i\right)\geq0.
\end{equation}

Suppose first that $\widetilde{\mathcal{L}}^2=0$. By Lemma \ref{lem:2.1}(1), the pencil defines a morphism $Y\to\P^1$ whose general fiber is a smooth rational curve. Exactly one component $E_j$ of $E$ is a section, and every other component is vertical. Consequently,
\[
 \widetilde{\mathcal{L}}\cdot\left(K_Y+\frac{\gamma_p}{a}\widetilde{\mathcal{L}}+\Theta\right)
 =-2+\gamma_pz_j\leq-2+\beta_p<0,
\]
contrary to \eqref{eq:2.12}.

Suppose next that $\widetilde{\mathcal{L}}^2>0$, and let $q$ be the base point of $\widetilde{\mathcal{L}}$. By Lemma \ref{lem:2.2} and \eqref{eq:2.12}, the pair $(Y,\frac{\gamma_p}{a}\widetilde{\mathcal{L}}+\Theta)$ is not log canonical at $q$. Since $E$ has simple normal crossings, at most two components of $\Theta$ pass through $q$, and each has coefficient at most $\beta_p<1$. Applying Lemma \ref{lem:2.3} at $q$, with a transverse smooth curve of coefficient zero if necessary, yields
\[
\widetilde{\mathcal{L}}^2>
\frac{4a^2}{\gamma_p^2}(1-\beta_p)^2.
\]
On the other hand, \eqref{eq:2.11} and \eqref{eq:2.10} imply
\[
 \widetilde{\mathcal{L}}^2\leq a^2h\leq\frac{4a^2}{\gamma_p^2}(1-\beta_p)^2,
\]
which is a contradiction. Therefore, $p$ cannot be the base point of a pencil associated to a cylinder.
\end{proof}

\section{Du Val del Pezzo surfaces}\label{sect:3}
For a Du Val del Pezzo surface $S$, the minimal resolution $\pi\colon Y\to S$ is crepant. Hence, $-K_Y=\pi^*(-K_S)$ is nef and big, and $Y$ is a weak del Pezzo surface. By \cite[Lemma 2.7]{MS24}, the surface $Y$, and therefore $S$, is rational. Moreover, Du Val surface singularities are $\Q$-factorial. In particular, Lemma \ref{lem:2.1} holds for $S$ when $\rho(S)=1$.

\subsection{The nonexistence of cylinders}
For a point of type ${\sf D}_4$, the exceptional divisor is not a chain. We use the following consequence of inversion of adjunction \cite[Theorem 5.1]{Tan18}. The notation for a minimal resolution and its anticanonical divisor is that of Subsection \ref{subsec:2.1}.

\begin{lemma}\label{lem:3.1}
    Let $T\subset V$ be a smooth curve through a smooth point $p$ of a surface $V$, and let $\Omega\geq 0$ be a $\Q$-divisor not containing $T$. If $(\Omega\cdot T)_p\leq 1$, then $(V,T+\Omega)$ is log canonical at $p$. 
\end{lemma}
\begin{proof}
    On a smooth neighborhood of $p$, adjunction gives $\mathrm{Diff}_T(\Omega)=\Omega|_T$, whose coefficient at $p$ is $(\Omega\cdot T)_p\leq 1$. Thus, $(T,\mathrm{Diff}_T(\Omega))$ is log canonical near $p$. The assertion then follows from \cite[Theorem 5.1]{Tan18}.
\end{proof}

\begin{remark}
    Lemma \ref{lem:3.1} also holds for $\Omega=c\mathcal{M}+B$, where $c\in \Q_{\geq 0}$, $\mathcal{M}$ is movable, and $B$ is an effective $\Q$-divisor not containing $T$. Indeed, fix a prime divisor $F$ over $V$ whose center contains $p$. Choose a general member $M\in\mathcal M$ such that $\ord_F(M)=\ord_F(\mathcal M)$ and the local intersection $(cM+B)\cdot T$ at $p$ equals its general value. Applying Lemma \ref{lem:3.1} to $cM+B$ gives $a_F(V,T+c\mathcal M+B)\geq-1$. Since $F$ was arbitrary, the assertion follows.
\end{remark}

\begin{lemma}\label{lem:3.3}
    Let $S$ be a Du Val del Pezzo surface with $(-K_S)^2=1$, and let $p\in S$ be a point of type ${\sf D}_4$. Let $\mathcal{L}$ be a movable linear system satisfying $\mathcal{L}\equiv -aK_S$ for some $a\in \Q_{>0}$. Then $\left(S,\frac{1}{a}\mathcal{L}\right)$ is log canonical at $p$.
\end{lemma}
\begin{proof}
    Recall that $H=\pi^{\ast}(-K_S)=-K_Y$. Label the exceptional curves over $p$ as $E_0,E_1,E_2,E_3$ so that $E_0$ meets each of $E_1,E_2,E_3$ and the latter three curves are pairwise disjoint. In the order $E_0,E_1,E_2,E_3$, the negative intersection matrix and the inverse matrix are
    \begin{equation}\label{eq:3.1}
    C=\begin{pmatrix}
    2&-1&-1&-1\\
    -1&2&0&0\\
    -1&0&2&0\\
    -1&0&0&2
    \end{pmatrix},\qquad
    G=C^{-1}=
    \begin{pmatrix}
    2&1&1&1\\   
    1&1&\frac{1}{2}&\frac{1}{2}\\
    1&\frac{1}{2}&1&\frac{1}{2}\\
    1&\frac{1}{2}&\frac{1}{2}&1
    \end{pmatrix}.
    \end{equation}

    We construct a smooth rational $(-1)$-curve $T\subset Y$ meeting $E_0$ once and disjoint from every other curve contracted by $\pi$. By Riemann--Roch, we have $\chi(Y,\OO_Y(H))=2$, and Serre duality gives $h^2(Y,\OO_Y(H))=h^0(Y,\OO_Y(-2H))=0$. Hence, $h^0(Y,\OO_Y(H))\geq 2$. Since $H=\pi^{\ast}(-K_S)$ and $\pi_{\ast}\OO_Y=\OO_S$, there is an effective Cartier divisor $T_S\in |-K_S|$ passing through $p$. Since $-K_S$ is an ample Cartier divisor, $(-K_S)\cdot C$ is a positive integer for every irreducible curve $C\subset S$. As $(-K_S)\cdot T_S=1$, the divisor $T_S$ is irreducible and reduced.

    Let $T$ be the strict transform of $T_S$. Write the exceptional part of $\pi^{\ast}T_S$ over $p$ as $\sum_{i=0}^3 m_iE_i$, and let $d_i=T\cdot E_i$. The vector $\bm{d}=(d_0,d_1,d_2,d_3)^{\rm t}$ is a nonzero vector with nonnegative integer entries. The equalities $\pi^{\ast}T_S\cdot E_i=0$ give $\bm{m}=G\bm{d}$, where $\bm{m}=(m_0,m_1,m_2,m_3)^{\rm t}$ is integral since $T_S$ is Cartier. For every other singular point $v$, write $\bm{m}_v$ for the coefficients of $\pi^{\ast}T_S$ over $v$ and $C_v$ for the negative intersection matrix. By orthogonality and adjunction, we have
    \[
    T^2=1-\bm{m}^{\rm t}C\bm{m}-\sum_{v\neq p} \bm{m}_v^{\rm t}C_v\bm{m}_v,\quad T^2=2p_a(T)-1\geq -1.
    \]
    Thus, we have $\bm{d}^{\rm t}G\bm{d}\leq 2$. Let $\bm{e}_0,\bm{e}_1,\bm{e}_2,\bm{e}_3$ be the coordinate vectors, with coordinates ordered by $E_0,E_1,E_2,E_3$. We determine the possible values of $\bm{d}$ from \eqref{eq:3.1}.

    Suppose first that $d_0\geq 1$. Since all entries of $G$ are nonnegative and $G_{00}=2$, we have $\bm{d}^{\rm t}G\bm{d}\geq 2d_0^2$. Thus, $d_0=1$. If $d_i\geq 1$ for some $i\in \{1,2,3\}$, then $\bm{d}^{\rm t}G\bm{d}\geq G_{00}+2G_{0i}+G_{ii}=2+2+1=5$, which is impossible. Hence, $\bm{d}=\bm{e}_0$.

    Suppose next that $d_0=0$. If two of $d_1,d_2,d_3$ are positive, say $d_i,d_j\geq 1$ with $i\neq j$, then $\bm{d}^{\rm t}G\bm{d}\geq G_{ii}+2G_{ij}+G_{jj}=1+1+1=3$, again a contradiction. If $d_i\geq 2$ for some $i\in \{1,2,3\}$, then $\bm{d}^{\rm t}G\bm{d}\geq 4G_{ii}=4$, which is also impossible. Therefore, in this case we have $\bm{d}=\bm{e}_i$ for some $i\in \{1,2,3\}$. 

    We have thus shown that $\bm{d}\in \{\bm{e}_0,\bm{e}_1,\bm{e}_2,\bm{e}_3\}$. On the other hand, $\bm{m}=G\bm{d}$ is an integral vector. For $i=1,2,3$, two entries of $G\bm{e}_i$ are equal to $\frac{1}{2}$. Hence, $G\bm{e}_i$ is not integral. Therefore, $\bm{d}=\bm{e}_0$, and $\bm{m}=G\bm{e}_0=(2,1,1,1)^{\rm t}$.

    All vectors $\bm{m}_v$ with $v\neq p$ vanish, and $T^2=-1$, $p_a(T)=0$. Consequently, $T$ is a smooth rational $(-1)$-curve, $H\cdot T=1$, $T\cdot E_0=1$, and $T$ is disjoint from the remaining exceptional curves.

    Let $\mathcal{B}\coloneqq \frac{1}{a}\widetilde{\mathcal{L}}$. For a general member $M\in \mathcal{L}$ with strict transform $\widetilde{M}$, write $\frac{1}{a}\pi^{\ast}M=\frac{1}{a}\widetilde{M}+\sum_{i=0}^3 x_iE_i$ near $\pi^{-1}(p)$. Set $\bm{x}\coloneqq (x_0,x_1,x_2,x_3)^{\rm t}$ and $b_i\coloneqq \mathcal{B}\cdot E_i\geq 0$. Then the equation $C\bm{x}=(b_0,b_1,b_2,b_3)^{\rm t}$ gives
    \[
    b_0=2x_0-x_1-x_2-x_3,\quad b_i=2x_i-x_0\quad (1\leq i\leq 3).
    \]
    Since every entry of $G$ is positive, $x_i\geq 0$. Moreover, since $\mathcal{B}$ has no fixed component and $T$ meets only $E_0$ among the exceptional curves, $0\leq \mathcal{B}\cdot T=1-x_0$. The inequalities $b_i\geq 0$ then give $\frac{x_0}{2}\leq x_i\leq x_0$ for $1\leq i\leq 3$.

    At a point of $E_0\setminus (E_1\cup E_2\cup E_3)$, the local intersection of $\mathcal{B}$ with $E_0$ is at most $b_0\leq \frac{x_0}{2}\leq 1$. At $E_0\cap E_i$, the local intersection number of $E_0$ with the remaining boundary is at most
    \[
    b_0+x_i=2x_0-\sum_{\substack{1\leq j\leq3\\j\neq i}}x_j\leq x_0\leq 1.
    \]
    After increasing the coefficient of $E_0$ to one, we conclude from Lemma \ref{lem:3.1} that the pair $(Y,\mathcal{B}+\sum_{i=0}^3 x_iE_i)$ is log canonical at every point of $E_0$. At a point of $E_i\setminus E_0$ for $i=1,2,3$, the local intersection with $E_i$ is at most $b_i=2x_i-x_0\leq x_0\leq 1$. After increasing the coefficient of $E_i$ to one, we again obtain log canonicity at every point of $E_i\setminus E_0$ by Lemma \ref{lem:3.1}. Since $\pi$ is crepant, $\left(S,\frac{1}{a}\mathcal{L}\right)$ is log canonical at $p$.
\end{proof}

\subsection{Curves on the minimal resolution}
Throughout this section, $\Sigma$ denotes the set of six Dynkin types listed in Theorem \ref{thm:1.2}. 

\begin{proposition}\label{prop:3.4}
    Let $S$ be a Du Val del Pezzo surface of Picard rank one over $\Bbbk$, and suppose that its Dynkin type occurs in characteristic zero. Let $S_{\C}$ be a Du Val del Pezzo surface of Picard rank one over $\C$ with the same Dynkin type. There is a bijection between the configurations of $(-1)$-curves and $(-2)$-curves on the minimal resolutions of $S$ and $S_{\C}$ that preserves their self-intersection and pairwise intersection numbers.
\end{proposition}
\begin{proof}
    Let $Y_{\C}\to S_{\C}$ be the minimal resolution. Recall that $h=H^2=(-K_S)^2=(-K_Y)^2$ is the degree of the weak del Pezzo surface $Y$. If $h=9$, then $Y\simeq\P^2_{\Bbbk}$ and $Y_{\C}\simeq\P^2_{\C}$. Neither surface contains a $(-1)$-curve or a $(-2)$-curve. We may therefore assume $h\leq8$. Except for $\P^1\times\P^1$ and $\F_2$, a weak del Pezzo surface is obtained from $\P^2$ by successive blow-ups; see \cite[Lemma 2.7]{MS24}. For $h<8$, its Picard group has a basis $\ell,e_1,\dots,e_{9-h}$ with
    \[
    \ell^2=1,\quad \ell\cdot e_i=0,\quad e_i\cdot e_j=-\delta_{ij},\quad K_Y=-3\ell+\sum_{i=1}^{9-h} e_i,
    \]
    where $\delta_{ij}$ is the Kronecker delta, $\ell$ is the pullback of the class of a line in $\P^2$, and $e_i$ is the class of the total transform of the exceptional curve of the $i$th blow-up. The Picard group of the minimal resolution $Y_{\C}$ of $S_{\C}$ has a basis with the same intersection numbers and canonical class.

    Every $(-2)$-curve on $Y$ has intersection zero with $H$, and is therefore contracted by $\pi$. Conversely, every component of $\Exc(\pi)$ is a $(-2)$-curve. Since $\rho(S)=1$, the number of these components is $9-h$, and their classes span $K_Y^{\perp}\otimes \R$.

    The case $h=7$ does not occur. Indeed, in this case $\rho(Y)=3$. Thus, the classes of the exceptional curves span the two-dimensional real vector space $K_Y^\perp\otimes_{\Z} \R$, whereas the root system of a weak del Pezzo surface of degree seven has type $\mathsf A_1$.

    If $h=8$, then $\rho(Y)=2$ and $\Exc(\pi)$ consists of one $(-2)$-curve. By \cite[Lemma 2.7]{MS24}, a weak del Pezzo surface of degree eight is isomorphic to $\P^1\times\P^1$, $\F_1$, or $\F_2$. The first two surfaces have no $(-2)$-curves. Hence, $Y\simeq\F_2$, and $\Exc(\pi)$ is the unique $(-2)$-section. The same argument holds for $Y_{\C}$, proving the assertion in degree eight.

    For $1\leq h\leq 6$, consider the root system
    \[
    \mathcal{R}_h\coloneqq \{v\in \Pic(Y)\mid v^2=-2, K_Y\cdot v=0\},
    \]
    with the positive definite form given by the negative of the intersection form. For $h=6,5,4,3,2,1$, the types of $\mathcal{R}_h$ are, respectively, ${\sf A}_2+{\sf A}_1$, ${\sf A}_4$, ${\sf D}_5$, ${\sf E}_6$, ${\sf E}_7$, and ${\sf E}_8$. The exceptional curves form a simple root basis for a root subsystem of full rank.

    Let $\mathcal S$ and $\mathcal S_{\C}$ be the root subsystems generated by the exceptional curves, after the two Picard lattices have been written in the bases above. Both subsystems have full rank. In the notation of \cite[Definition 2]{Lub14}, this gives $C_{\mathcal R_h}(\mathcal S)=C_{\mathcal R_h}(\mathcal S_{\C})=\emptyset$ and $K_{\mathcal R_h}(\mathcal S)=K_{\mathcal R_h}(\mathcal S_{\C})=\mathcal R_h$. Thus, the two additional invariants in \cite[Proposition 3]{Lub14} agree. Since the Dynkin types also agree, that proposition gives Weyl conjugacy when $\mathcal R_h$ is irreducible. For $\mathcal R_h$ of type ${\sf A}_2+{\sf A}_1$, the only root subsystem of full rank is $\mathcal R_h$ itself. Composing with an element of the Weyl group of $\mathcal S_{\C}$ then sends the chosen simple root basis of $\mathcal S$ to that of $\mathcal S_{\C}$. Each reflection in a root $r$ acts on the Picard group by $v\mapsto v+(v\cdot r)r$. Consequently, there is an integral isometry
    \[
    \Phi\colon \Pic(Y)\longrightarrow \Pic(Y_{\C})
    \]
    which preserves canonical classes and sends the classes of the exceptional curves to the classes of the exceptional curves of $Y_{\C}$. 

    We next show that $\Phi$ maps the classes of $(-1)$-curves on $Y$ to those on $Y_{\C}$. A class $v\in \Pic(Y)$ is represented by a $(-1)$-curve if and only if
    \begin{equation}\label{eq:3.2}
        v^2=-1,\quad H\cdot v=1,\quad v\cdot T\geq 0 \text{ for every irreducible component } T\subset \Exc(\pi).
    \end{equation}
    By the Riemann--Roch formula, we have $\chi(Y,\OO_Y(v))=1$. Since $H\cdot (K_Y-v)=-h-1<0$ and $H$ is nef, Serre duality gives $h^2(Y,\OO_Y(v))=0$. Hence, there exists an effective divisor $D\sim v$. The equality $H\cdot D=1$ implies that $D=C+Z$, where $C$ is an irreducible component appearing with multiplicity one and satisfying $H\cdot C=1$, and $Z\geq 0$ is supported on $\Exc(\pi)$. By adjunction, we have $C^2=2p_a(C)-1\geq -1$. Moreover, \eqref{eq:3.2} gives $D\cdot Z\geq 0$. If $Z\neq 0$, then negative definiteness of the exceptional intersection matrix gives
    \[
    -1=D^2=C^2+2D\cdot Z-Z^2>-1,
    \]
    which is a contradiction. Thus, $D=C$, and adjunction gives $p_a(C)=0$. The curve $C$ is a smooth rational $(-1)$-curve, and is the unique irreducible curve in its class, since two distinct irreducible curves cannot have negative intersection. Hence, we obtain \eqref{eq:3.2}. The isometry $\Phi$ preserves all three conditions in \eqref{eq:3.2}. Thus, $\Phi$ induces the required bijection and preserves every pairwise intersection number.
\end{proof}

\subsection{Cylinders on surfaces of Picard rank one}

\begin{proposition}\label{prop:3.5}
Let $S$ be a Du Val del Pezzo surface of Picard rank one over $\Bbbk$. If $\mathrm{Dyn}(S)\notin\Sigma$, then $S$ contains a cylinder.
\end{proposition}

\begin{proof}
Let $E\coloneqq\Exc(\pi)_{\rm red}$. If $h=9$, then $S\simeq\P^2$ contains $\A^2$. If $h=8$, then $Y\simeq\F_2$, and $Y\setminus \Supp(E+F)\simeq \A^2$ for a fiber $F$ of $\F_2\to \P^1$. The case $h=7$ does not occur; see \cite[Theorem 1.1 and Table 1]{KN23}. The remaining degrees are $1\leq h\leq 6$.

By \cite[Theorem 1.1 and Remark 1.2]{KN23}, the Dynkin type $\mathrm{Dyn}(S)$ occurs in characteristic zero. Proposition \ref{prop:3.4} identifies the configurations of $(-1)$-curves and $(-2)$-curves on $Y$ with those in \cite[Tables 3.1, 5.1, and 5.5]{Saw24}, including their intersection numbers. We first treat the types other than $\mathsf E_6+\mathsf A_2$, $\mathsf A_5+\mathsf A_2+\mathsf A_1$, and $2\mathsf A_4$ by the contraction sequences specified below.

For $h\geq3$, we use \cite[Table 3.1 and Section 3.2]{Saw24}. For $h=2$, we use \cite[proof of Lemma 5.4]{Saw24} for $\mathsf D_4+3\mathsf A_1$, \cite[proof of Lemma 5.5]{Saw24} for $\mathsf A_7$, $\mathsf D_6+\mathsf A_1$, and $\mathsf E_7$, and \cite[proof of Lemma 5.22(1)]{Saw24} for $\mathsf A_5+\mathsf A_2$ and $2\mathsf A_3+\mathsf A_1$. For $h=1$, we use \cite[proof of Lemma 5.5]{Saw24} for $\mathsf A_8$, $\mathsf D_8$, $\mathsf E_7+\mathsf A_1$, and $\mathsf E_8$, and \cite[proof of Lemma 5.22(2)]{Saw24} for $\mathsf A_7+\mathsf A_1$. For $\mathsf D_6+2\mathsf A_1$, the first contraction in \cite[Table 5.1]{Saw24} gives a weak del Pezzo surface of degree two and type $\mathsf D_4+3\mathsf A_1$; we then use the construction in \cite[proof of Lemma 5.4]{Saw24}. For type $\mathsf D_5+\mathsf A_3$, we reduce to a weak del Pezzo surface of degree two below. These are all the remaining types by \cite[Theorem 1.1 and Table 1]{KN23}.

Let $\sigma$ be the contraction of a smooth rational $(-1)$-curve $J$, and let $D,D'$ be irreducible curves distinct from $J$. Then
\begin{equation}\label{eq:3.3}
\begin{aligned}
(\sigma_*D)^2&=D^2+(D\cdot J)^2,\\
\sigma_*D\cdot\sigma_*D'&=D\cdot D'+(D\cdot J)(D'\cdot J),\\
p_a(\sigma_*D)&=p_a(D)+\frac{(D\cdot J)(D\cdot J-1)}{2}.
\end{aligned}
\end{equation}

We use the full contraction sequences described above, including those on intermediate weak del Pezzo surfaces. On the corresponding minimal resolution in characteristic zero, the strict transform of every curve contracted in these sequences has negative self-intersection. By nefness of the anticanonical divisor and adjunction, it is a $(-1)$-curve or a $(-2)$-curve. Proposition \ref{prop:3.4} provides the corresponding curves on $Y$. We apply \eqref{eq:3.3} inductively. At each step, the next curve to be contracted is integral, and its self-intersection number and arithmetic genus agree with those in the corresponding characteristic-zero construction. These values are $-1$ and $0$, respectively. The curve is therefore smooth and rational, and Castelnuovo's contraction theorem applies over $\Bbbk$. The same formulas determine the intersection numbers and arithmetic genera of the remaining curves. Thus, these contraction sequences can be carried out over $\Bbbk$.

For $h=1$ and $\mathrm{Dyn}(S)=\mathsf D_5+\mathsf A_3$, contract the $(-1)$-curve $\widetilde E$ marked in \cite[Table 5.5]{Saw24}. Let $\sigma\colon Y\to\overline Y$ denote this contraction. For every irreducible curve $\overline D\subset\overline Y$ with strict transform $D$, we have $(-K_{\overline Y})\cdot\overline D=(-K_Y)\cdot D+\widetilde E\cdot D\geq0$. Moreover, $K_{\overline Y}^2=K_Y^2+1=2$. Thus, $\overline Y$ is a weak del Pezzo surface of degree two. Formula \eqref{eq:3.3} gives the configuration of type $2\mathsf A_3+\mathsf A_1$. The transforms of the curves in the original configuration give the corresponding degree-two configuration in the same table. We then use the contraction sequence in \cite[proof of Lemma 5.22(1)]{Saw24}.

For $\mathrm{Dyn}(S)=\mathsf E_6+\mathsf A_2$, we give a contraction sequence directly. Using the Weyl conjugacy in the proof of Proposition \ref{prop:3.4}, choose a lattice marking $\ell,e_1,\dots,e_8$ with $\ell^2=1$, $e_i^2=-1$, all other pairwise intersections zero, and $K_Y=-3\ell+\sum_i e_i$. With a suitable labeling, the exceptional curves $E_1,\dots,E_8$ have respective classes
\[
e_1-e_2,\ e_2-e_3,\ e_3-e_4,\ e_4-e_5,\ e_5-e_6,\ e_7-e_8,
\ \ell-e_1-e_2-e_3,\ -3\ell+e_1+\cdots+e_7+2e_8.
\]
Here the marking is an isometry of Picard lattices, not a choice of a morphism to $\P^2$. The classes $\ell-e_7-e_8$ and $3\ell-e_1-\cdots-e_5-2e_7-e_8$ satisfy \eqref{eq:3.2}. Hence, they are represented by $(-1)$-curves $C_1,C_2$. Contract successively the transforms of $C_1,E_7,E_3,E_2,E_1,C_2,E_5$. Formula \eqref{eq:3.3} shows that each is a smooth rational $(-1)$-curve when contracted. The resulting surface is smooth and rational of Picard rank two, and the image of $E_8$ is a $(-2)$-curve. Hence, it is isomorphic to $\F_2$. The images of $E_8,E_4,E_6$ are its negative section, a section disjoint from it, and a fiber, respectively.

For these contraction sequences, write $\tau\colon Y\to V$ for the resulting morphism. The configurations in \cite[Tables 3.1, 5.1, and 5.5]{Saw24}, together with \eqref{eq:3.3}, show that the one-dimensional part of $\tau(E)$ consists of one or two lines on $\P^2$, a negative section and fibers on $\F_2$, two disjoint sections and at most one fiber on $\F_2$, or two sections meeting once on $\F_2$ with no additional fiber. In the last case, the sections have classes $M$ and $M+3F$, where $M$ is the negative section and $F$ is a fiber. Their intersection number is one. Hence, they meet transversely at a single point. In the case of two disjoint sections their classes are $M$ and $M+2F$. The explicit construction for $\mathsf E_6+\mathsf A_2$ has the same boundary form.

The complements of one and two lines in $\P^2$ are $\A^2$ and $\A^1\times\mathbb G_m$, respectively. When the image contains a single section on $\F_2$, we remove this section and all fibers appearing in the image, removing one fiber if no fiber appears. The complement is an affine-line bundle over an affine open subset of $\P^1$. Both the associated line bundle and the torsor class are trivial. Hence, it is a cylinder. When the image contains two disjoint sections, we remove the fiber appearing in the image of $E$, or an arbitrary fiber if none appears. The base is then $\A^1$, and the complement of the two sections is $\A^1\times\mathbb G_m$. When the two sections meet once, no additional fiber occurs in the image of $E$. Removing the fiber through their intersection again gives $\A^1\times\mathbb G_m$.

For this argument, it suffices to choose a cylinder avoiding the one-dimensional components of the image of the exceptional locus. More precisely, let $g\colon W\to Y$ and $\tau\colon W\to V$ be birational morphisms, and suppose that $V$ has a cylinder $U\simeq\A^1\times Z$ disjoint from every curve in $\tau(\Exc(\pi\circ g))$. The images in $V$ of $\Exc(\tau)$ and of the components of $\Exc(\pi\circ g)$ contracted by $\tau$ form a finite set. Remove from $Z$ the images under the projection $U\to Z$ of the points of this finite set lying in $U$. The resulting smaller cylinder lifts isomorphically to an open subset of $W$ disjoint from $\Exc(\pi\circ g)$, and hence maps isomorphically onto an open subset of $S$. Applying this with $W=Y$ proves the assertion in all the cases above.

Suppose next that $\mathrm{Dyn}(S)=\mathsf A_5+\mathsf A_2+\mathsf A_1$. Write the exceptional configuration as the disjoint union of the chains $E_1-E_2$ and $E_3-E_4-E_5-E_6-E_8$ and the curve $E_7$. By Proposition \ref{prop:3.4} and \cite[Table 5.5]{Saw24}, there are $(-1)$-curves $C_1,C_2$ with $C_1\cdot E_1=C_1\cdot E_4=C_2\cdot E_5=C_2\cdot E_7=1$ and all their other intersections with exceptional curves equal to zero, and $C_1\cdot C_2=0$. Neither curve passes through $E_4\cap E_5$. Let $g\colon W\to Y$ be the blow-up of this transverse intersection, with exceptional curve $G$. Contract successively the strict transforms of $C_1,E_1,E_2,C_2,E_7,E_5,E_6,E_8$. Formula \eqref{eq:3.3} shows that each is a smooth rational $(-1)$-curve when contracted. The resulting surface is $\F_2$, and the images of $E_3,E_4,G$ have intersection matrix 
$\left(
\begin{smallmatrix}
-2&1&0\\1&0&1\\0&1&2
\end{smallmatrix}
\right)$. These images are the negative section, a fiber, and a section disjoint from the negative section. They are the only one-dimensional components of the image of $\Exc(\pi\circ g)$. Their complement is $\A^1\times\mathbb G_m$. By the preceding argument, a smaller cylinder descends to $S$.

Finally, suppose that $\mathrm{Dyn}(S)=2\mathsf A_4$. Label the exceptional chains as $A_1-A_2-A_3-A_4$ and $B_1-B_2-B_3-B_4$. By Proposition \ref{prop:3.4} and \cite[Table 5.5]{Saw24}, there is a $(-1)$-curve $C$ meeting $A_2$ and $B_4$ transversely once and disjoint from all the other exceptional curves. Blow up successively the transverse intersections $A_2\cap A_3$, $A_2\cap F_1$, $A_2\cap F_2$, $A_2\cap F_3$, and $F_3\cap F_4$. Here $F_i$ denotes the exceptional curve of the $i$th blow-up, and we use the same notation for the strict transforms of the curves. Each blow-up is centered at the intersection of exactly two components of a simple normal crossings divisor, and none of these centers lies on $C$. Let $g\colon W\to Y$ be the resulting morphism, and put
\[
D\coloneqq A_1+2A_2+2B_1+4B_2+6B_3+8B_4+10C+F_4.
\]
On $W$, we have $A_2^2=-6$, $C^2=-1$, and all other components of $D$ have self-intersection $-2$. Their intersection graph is the chain $B_1-B_2-B_3-B_4-C-A_2-A_1$, with the additional edge $A_2-F_4$. A direct computation of the intersection numbers gives $D\cdot T=0$ for every component $T$ of $D$. Since $D$ is effective, it is nef and $D^2=0$. Adjunction gives $K_W\cdot D=-2$. Moreover, $D\cdot F_5=1$, and $D$ has intersection zero with all the other curves contracted by $\pi\circ g$.

We show that $|D|$ is a base-point-free pencil. Since $D$ is nef and $D\cdot(K_W-D)=-2$, the divisor $K_W-D$ is not effective. Hence, $h^0(W,\OO_W(K_W-D))=0$. Since $W$ is rational, Riemann--Roch and Serre duality yield $h^0(W,\OO_W(D))\geq2$. Write $|D|=|M|+Z$, where $Z$ is the fixed part and $|M|$ is the movable part. Since $D$ is nef and $D^2=0$, we have $D\cdot M=D\cdot Z=0$. Moreover, $M^2\geq0$ because $|M|$ has no fixed components. The Hodge index theorem gives $M\equiv aD$. Intersecting with an ample divisor and using $M\neq0$ and $Z\geq0$, we obtain $0<a\leq1$. Since $M\cdot F_5=a$ is a positive integer, $a=1$. Hence, $Z\equiv0$, and effectivity implies $Z=0$. Thus, $|D|$ has no fixed component, and $D^2=0$ implies that it has no base point.

Since $D^2=0$ and $h^0(W,\OO_W(D))\geq2$, the morphism defined by $\lvert D\rvert$ has one-dimensional image. Let $f\colon W\to B$ be the Stein factorization of the morphism defined by $|D|$. Then $D\sim f^*A$ for an ample divisor $A$ on the smooth projective curve $B$. The equality $D\cdot F_5=1$ implies $\deg A=1$ and $\deg(F_5\to B)=1$. Therefore, $B\simeq\P^1$, $|D|$ is a pencil, and $F_5$ is a section. The generic fiber is a regular curve with arithmetic genus zero by adjunction; the section provides a $\Bbbk(B)$-rational point on the generic fiber. Hence, this fiber is $\P^1_{\Bbbk(B)}$.

Every curve contracted by $\pi\circ g$ other than $F_5$ is vertical. Over a nonempty affine open subset $B^\circ\subset B$ avoiding their images, the isomorphism on the generic fiber extends to $f^{-1}(B^\circ)\simeq\P^1\times B^\circ$, with $F_5$ corresponding to $\{\infty\}\times B^\circ$. Removing this section gives a cylinder $\A^1\times B^\circ$ disjoint from $\Exc(\pi\circ g)$. It descends isomorphically to a cylinder in $S$.
\end{proof}

\subsection{Blowing up a surface of type $7{\sf A}_1$}
\begin{lemma}\label{lem:3.6}
    Let $S$ be a Du Val del Pezzo surface of Picard rank one and Dynkin type $7{\sf A}_1$, and let $q\in S$ be a smooth point. Let $\beta\colon S_q\to S$ be the blow-up of $q$. Then $-K_{S_q}$ is nef and big and is not ample. Let $\pi_q\colon Y_q\to S_q$ be the minimal resolution. Then $Y_q$ contains a $(-2)$-curve not contracted by $\pi_q$.
\end{lemma}
\begin{proof}
    We first note that the surface has degree two and occurs only in characteristic two. By \cite[Proposition 3.14]{KN23}, the defining equation is
    \[
    S=\{w^2+xyz(x+y+z)=0\}\subset \P(1,1,1,2).
    \]
    The anticanonical morphism $\varphi\colon S\to \P^2$, given by $[x:y:z:w]\mapsto [x:y:z]$, is finite and purely inseparable of degree two. Let $A$ be the exceptional curve of $\beta$.

    We prove that $-K_{S_q}=\beta^*(-K_S)-A$ is nef. Let $C\subset S$ be an irreducible curve through $q$. Since $\varphi$ is finite, there is a line $\ell\subset\P^2$ through $\varphi(q)$ which does not contain $\varphi(C)$. The divisor $D_{\ell}\coloneqq\varphi^*\ell\in|-K_S|$ contains $q$ and does not contain $C$. Since $q$ is smooth, local intersection multiplicities give
\[
(-K_S)\cdot C=D_{\ell}\cdot C
\geq(D_{\ell}\cdot C)_q\geq\mult_q C.
\]
Thus, the strict transform $C_q$ satisfies $(-K_{S_q})\cdot C_q\geq0$. Curves not passing through $q$ have positive anticanonical degree after the blow-up, and $(-K_{S_q})\cdot A=1$. Hence, $-K_{S_q}$ is nef. The equality $(-K_{S_q})^2=(-K_S)^2-1=1$ proves bigness. For the minimal resolution $\pi_q\colon Y_q\to S_q$, crepancy gives $-K_{Y_q}=\pi_q^*(-K_{S_q})$, which is also nef and big.

    In an affine chart containing $q$, write the equation as $w^2=f(u,v)$. Smoothness at $q$ gives $(f_u(q),f_v(q))\neq (0,0)$, and the differential of $\varphi$ at $q$ has rank one. Thus, there is a line through $\varphi(q)$ whose local equation pulls back to an element of $\mathfrak{m}_q^2$. Let $D\in |-K_S|$ be its pullback and let $m\coloneqq \mult_q D\geq 2$. The strict transform $D_q$ is an effective divisor, and 
    \[
    (-K_{S_q})\cdot D_q=(\beta^*(-K_S)-A)\cdot (\beta^*D-mA)=2-m.
    \]
    Since $-K_{S_q}$ is nef, we obtain that $m=2$ and $(-K_{S_q})\cdot D_q=0$. Hence, $-K_{S_q}$ is not ample.

    Let $R\subset Y_q$ be the strict transform of an arbitrary irreducible component of $D_q$. Since $D_q$ is effective, $-K_{S_q}$ is nef, and $(-K_{S_q})\cdot D_q=0$, the projection formula gives $(-K_{Y_q})\cdot R=0$. Moreover, since $(-K_{Y_q})^2=1$, the Hodge index theorem gives $R^2<0$. By adjunction, we have $R^2=2p_a(R)-2$, and hence $R^2=-2$ and $p_a(R)=0$. Therefore, $R$ is a smooth rational $(-2)$-curve, and its image on $S_q$ is a curve. This proves the last assertion.
\end{proof}

\subsection{Birational contractions}
\begin{lemma}[\protect{\cite[Lemma 9(b)]{DT19}}]\label{lem:3.7}
    Let $\theta\colon X\to T$ be a birational morphism of normal projective surfaces. If $T$ contains a cylinder, then $X$ contains a cylinder.
\end{lemma}
\begin{proof}
Let $U\simeq\A^1\times Z\subset T$ be a cylinder, and let $Q\subset T$ be the finite set over which $\theta$ is not an isomorphism. Remove from $Z$ the finite set $\operatorname{pr}_Z(Q\cap U)$, and denote the resulting nonempty affine curve by $Z^{\circ}$. Then $\theta^{-1}(\A^1\times Z^{\circ})\simeq\A^1\times Z^{\circ}$ is a cylinder in $X$.
\end{proof}

\begin{lemma}\label{lem:3.8}
    Every birational contraction in a $K$-minimal model program on a Du Val del Pezzo surface gives another Du Val del Pezzo surface and strictly increases the anticanonical degree. Furthermore, if a birational morphism $\theta\colon X\to T$ between Du Val del Pezzo surfaces satisfies $K_T^2-K_X^2=1$, then $\theta$ is the blow-up of a smooth point of $T$.
\end{lemma}
\begin{proof}
    Let $\varphi\colon S\to T$ be a birational extremal contraction and let $C$ be its exceptional curve. The surface $T$ is $\Q$-factorial by \cite[Theorem 1.1]{Tan18}. Write
    \[
    K_S=\varphi^*K_T+aC,\quad a=\frac{K_S\cdot C}{C^2}>0.
    \]
    For the minimal resolution $\mu\colon V\to S$, we have $K_V=(\varphi\circ\mu)^*K_T+a\mu^*C$. The divisor $\mu^*C$ is effective and supported on the exceptional locus of $\varphi\circ\mu$. Thus, every discrepancy of this resolution of $T$ is nonnegative, and $T$ has Du Val singularities. Moreover,
    \begin{equation}\label{eq:3.4}
        K_T^2=K_S^2-a^2C^2>K_S^2.
    \end{equation}
    For an irreducible curve $D\subset T$, let $D_S$ be its strict transform. Then $(-K_T)\cdot D=(-K_S)\cdot D_S+a(C\cdot D_S)>0$. By the Nakai--Moishezon criterion and \eqref{eq:3.4}, $-K_T$ is ample. 

    Let $\pi_X\colon Y_X\to X$ and $\pi_T\colon Y_T\to T$ be the minimal resolutions. The resolution $\theta\circ\pi_X$ factors through $Y_T$, giving a birational morphism $\sigma\colon Y_X\to Y_T$. Since $\pi_X$ and $\pi_T$ are crepant, we have $K_{Y_T}^2-K_{Y_X}^2=1$. A birational morphism between smooth surfaces is a sequence of point blow-ups, each decreasing $K^2$ by one. Therefore, $\sigma$ is the blow-up of a point $q\in Y_T$.

    If $q$ belongs to a component $D$ of $\Exc(\pi_T)$, then its strict transform satisfies $(-K_{Y_X})\cdot \widetilde{D}=-1$, since $D^2=-2$ and $(-K_{Y_T})\cdot D=0$. This contradicts the nefness of $-K_{Y_X}=\pi_X^*(-K_X)$. Hence, $q\not\in \Exc(\pi_T)$ and $\pi_T(q)$ is smooth. The strict transforms on $Y_X$ of the exceptional curves of $\pi_T$ have anticanonical degree zero and are therefore contracted by $\pi_X$. The new exceptional curve has anticanonical degree one and is not contracted by $\pi_X$. Since $\theta\circ\pi_X=\pi_T\circ\sigma$, the curves contracted by $\pi_X$ are precisely the strict transforms of the exceptional curves of $\pi_T$. Contracting the exceptional curves shows that $\theta$ is the blow-up of $\pi_T(q)$.
\end{proof}

\begin{lemma}\label{lem:3.9}
    Let $f\colon S\to B$ be a Mori fiber space from a Du Val surface to a smooth projective curve. Then $S$ contains a cylinder.
\end{lemma}
\begin{proof}
    Let $K=\Bbbk(B)$ and let $F$ be the generic fiber of $f$. The curve $F$ is regular, since $S$ has only isolated singularities, and $H^0(F,\OO_F)=K$ since $f_*\OO_S=\OO_B$. By adjunction and the relative ampleness of $-K_S$, we have $\deg \omega_F<0$, hence $p_a(F)=0$. The anticanonical linear system embeds $F$ as a plane conic. By Tsen's theorem \cite[Theorem 59.67.10, Tag 03RD]{Stacks}, the conic has a $K$-rational point. A regular conic with a $K$-rational point is isomorphic to $\P^1_K$. Fix such an isomorphism. The isomorphism and its inverse extend over a nonempty open subset $B^{\circ}\subset B$; after shrinking $B^{\circ}$, the two compositions are the identity. Thus, $S\times_B B^{\circ}\simeq\P^1\times B^{\circ}$. Shrink $B^{\circ}$ further to make it affine. Removing $\{\infty\}\times B^{\circ}$ gives a cylinder $\A^1\times B^{\circ}\subset S$.
\end{proof}

\subsection{Proof of Theorem \ref{thm:1.2}}
\begin{proof}[Proof of Theorem \ref{thm:1.2}]
    Assume first that $\rho(X)=1$ and that $\mathrm{Dyn}(X)\in\Sigma$.
    Let $\pi\colon Y\to X$ be the minimal resolution. Since $Y$ is a smooth projective rational surface, we have $K_Y^2+\rho(Y)=10$. Moreover, since $\pi$ is crepant and $\rho(X)=1$, the degree $(-K_X)^2$ is one for $4\mathsf A_2$, $2\mathsf D_4$, $2\mathsf A_1+2\mathsf A_3$, $4\mathsf A_1+\mathsf D_4$, and $8\mathsf A_1$, and is two for $7\mathsf A_1$.

    Suppose that $X$ contains a cylinder. Lemma \ref{lem:2.1} gives an associated pencil $\mathcal{L}\equiv-aK_X$ with a unique base point $P$ at which $\left(X,\frac{1}{a}\mathcal{L}\right)$ is not log canonical. Corollary \ref{cor:2.6} excludes a smooth point $P$, since the degree is at most two. For a degree one surface, Corollary \ref{cor:2.6} excludes a point of type $\mathsf A_i$ with $1\leq i\leq3$, since $d([2^i])=i+1\leq4$. For a degree two surface of type $7\mathsf A_1$, Corollary \ref{cor:2.6} would require $2>\frac{4}{2}=2$ at a point of type $\mathsf A_1$, a contradiction. A point of type $\mathsf D_4$ in the stated list occurs on a degree one surface and is excluded by Lemma \ref{lem:3.3}. This excludes all possible base points. Hence, $X$ contains no cylinder.

    Conversely, suppose that $X$ contains no cylinder. Run a $K_X$-minimal model program as in \cite[Theorem 1.1]{Tan18}, and let $\theta\colon X\to X'$ be the composition of its birational contractions. By Lemma \ref{lem:3.8}, each intermediate surface is a Du Val del Pezzo surface. The canonical divisor of every intermediate surface is not nef. Therefore, the $K$-minimal model program terminates with a Mori fiber space. Lemma \ref{lem:3.7} shows that $X'$ contains no cylinder. When the base of the Mori fiber space is a curve, Lemma \ref{lem:3.9} gives a cylinder in $X'$. Hence, the base is a point and $\rho(X')=1$.

    Proposition \ref{prop:3.5} now gives $\mathrm{Dyn}(X')\in\Sigma$. The five types in $\Sigma$ other than $7\mathsf A_1$ have degree one, whereas $7\mathsf A_1$ has degree two. Since $X$ is a Du Val del Pezzo surface, $K_X^2$ is a positive integer. Thus, if $K_{X'}^2=1$, then $K_X^2=1$. The strict increase of the anticanonical degree under every nontrivial birational contraction, proved in Lemma \ref{lem:3.8}, shows that $\theta$ is an isomorphism.

    Suppose that $K_{X'}^2=2$. Then $X'$ has type $7\mathsf A_1$. Suppose that $\theta$ is nontrivial. Then $K_X^2=1$, and Lemma \ref{lem:3.8} shows that $\theta$ is the blow-up of a smooth point of $X'$. Lemma \ref{lem:3.6} then implies that $-K_X$ is not ample, a contradiction. Therefore, $\theta$ is an isomorphism in this case as well. We obtain $\rho(X)=1$ and $\mathrm{Dyn}(X)\in\Sigma$, which proves the converse.
\end{proof}

\section{Proof of the Main Theorem}\label{sect:4}
\subsection{Non-log-liftable surfaces}
Throughout this section, $\mathrm{char}(\Bbbk)=5$. By \cite[Theorem 1.2 and Table 1]{Nag25}, the log del Pezzo surfaces of Picard rank one which are not log liftable are listed in Table \ref{tab:1}. The anticanonical degrees in the last column can be computed from \eqref{eq:4.1} for rows I--III and Nos. 1--3 (see Subsection \ref{subsec:4.2}), and from Lemmas \ref{lem:4.2} and \ref{lem:4.3} for Nos. 4--21. For Nos. 5--21, the parameters $r,u$ and the type of $B$ are specified in Subsection \ref{subsec:4.3} and Table \ref{tab:2}. These surfaces are rational, as follows from the constructions in \cite[Example 4.2 and Lemmas 6.2--6.3]{Nag25}.

\renewcommand{\arraystretch}{1.4}
\setlength{\tabcolsep}{5pt}
\captionsetup{width=\textwidth}

\begin{xltabular}{\textwidth}{
|>{\centering\arraybackslash}p{0.06\textwidth}
|>{\raggedright\arraybackslash\small}X
|>{\centering\arraybackslash\footnotesize}p{0.28\textwidth}|
}

\caption{Log del Pezzo surfaces of Picard rank one in characteristic five which are not log liftable.}
\label{tab:1}\\
\hline
\textbf{No.}
& \textbf{Exceptional configuration} ($n\geq 0$, $m\geq 1$)
& \textbf{$(-K_S)^2$} \\
\hline
\endfirsthead

\hline
\textbf{No.}
& \textbf{Exceptional configuration} ($n\geq 0$, $m\geq 1$)
& \textbf{$(-K_S)^2$} \\
\hline
\endhead

\hline
\endfoot

I
& $2[3,2]+[3]+[2]+[2^4]$
& $\frac{2}{15}$ \\
\hline

II
& $[4,2]+[3,2^5]+[3,2]+[2]$
& $\frac{2}{455}$ \\
\hline

III
& $[2,4]+[2,3,2^2]+[3]+[2^4]$
& $\frac{5}{231}$ \\
\hline

1.
& $2[2^4]+[3]$
& $\frac{1}{3}$ \\
\hline

2.
& $2[2^4]+[2,4]$
& $\frac{1}{7}$ \\
\hline

3.
& $2[2^4]+[2]+[3]+[5]$
& $\frac{2}{15}$ \\
\hline

4.
& $2[2^4]+[2^n]+[2+n;[2],[3],[5]]$
& $\frac{n+1}{30n+29}$ \\
\hline

5.
& $2[2^4]+[2^n,3]+[3,2+n,5]$
& $\frac{(n+2)^2}{(2n+3)(15n+22)}$ \\
\hline

6.
& $2[2^4]+[2^n,4]+[2,2+n,5]$
& $\frac{(n+2)^2}{(3n+4)(10n+13)}$ \\
\hline

7.
& $2[2^4]+[2^n,6]+[2,2+n,3]$
& $\frac{(n+2)^2}{(5n+6)(6n+7)}$ \\
\hline

8.
& $2[2^4]+[4]+[2;[2],[3],[5]]$
& $\frac{1}{29}$ \\
\hline

9.
& $2[2^4]+[3,2^{m-1},3]+[2+m;[2],[3],[5]]$
& $\frac{m+1}{30m+29}$ \\
\hline

10.
& $2[2^4]+[4+\ell]+[2;[2],[2^{\ell}],[5]] \quad (1\leq \ell\leq 2)$
& $\frac{4}{(\ell+4)(3\ell+13)}$ \\
\hline

11.
& $2[2^4]+[2+\ell,2^{m-1},4]+[2+m;[2],[2^{\ell}],[5]]
   \quad (1\leq \ell\leq 2)$
& $\frac{((\ell+1)m+2)^2}
        {(10(\ell+1)m+3\ell+13)(3(\ell+1)m+\ell+4)}$ \\
\hline

12.
& $2[2^4]+[2,5]+[2;[2],[3],[5]]$
& $\frac{1}{29}$ \\
\hline

13.
& $2[2^4]+[2,3,2^{m-1},4]+[2+m;[2],[3],[5]]$
& $\frac{m+1}{30m+29}$ \\
\hline

14.
& $2[2^4]+[6+\ell]+[2;[2],[3],[2^{\ell}]] \quad (1\leq \ell\leq 4)$
& $\frac{4}{(\ell+6)(\ell+7)}$ \\
\hline

15.
& $2[2^4]+[2+\ell,2^{m-1},6]+[2+m;[2],[3],[2^{\ell}]]
   \quad (1\leq \ell\leq 4)$
& $\frac{((\ell+1)m+2)^2}
        {(5(\ell+1)m+\ell+6)(6(\ell+1)m+\ell+7)}$ \\
\hline

16.
& $2[2^4]+[2^{\ell},7]+[2;[2],[3],[2+\ell]] \quad (1\leq \ell\leq 3)$
& $\frac{(\ell+2)^2}{(6\ell+7)(7\ell+8)}$ \\
\hline

17.
& $2[2^4]+[2^{\ell},3,2^{m-1},6]+[2+m;[2],[3],[2+\ell]]
   \quad (1\leq \ell\leq 3)$
& $\frac{(\ell+2)^2(m+1)^2}
        {(5(\ell+2)m+6\ell+7)(6(\ell+2)m+7\ell+8)}$ \\
\hline

18.
& $2[2^4]+[3,7]+[2;[2],[3],[3,2]]$
& $\frac{4}{115}$ \\
\hline

19.
& $2[2^4]+[3,3,2^{m-1},6]+[2+m;[2],[3],[3,2]]$
& $\frac{5m+4}{5(30m+23)}$ \\
\hline

20.
& $2[2^4]+[2,8]+[2;[2],[3],[2,3]]$
& $\frac{3}{85}$ \\
\hline

21.
& $2[2^4]+[2,4,2^{m-1},6]+[2+m;[2],[3],[2,3]]$
& $\frac{5m+3}{5(30m+17)}$ \\
\hline

\end{xltabular}

\subsection{Local intersection inequality}\label{subsec:4.2}
For each singular point of a surface appearing in rows I--III or Nos. 1--3, the exceptional divisor on the minimal resolution is a chain. We compute the anticanonical degrees using the intersection matrices, and then apply Corollary \ref{cor:2.6}.

Let $\pi\colon Y\to S$ be the minimal resolution of a surface $S$ in rows I--III or Nos. 1--3. For each singular point indexed by $\nu$, let $C_{\nu}$ be the negative of the exceptional intersection matrix, let $s_{\nu}$ be the number of exceptional components, and let $\bm{v}_{\nu}\coloneqq \left((C_{\nu})_{ii}-2\right)_{1\leq i\leq s_{\nu}}^{\rm t}$. The equality $K_Y^2+\rho(Y)=10$ on $Y$ gives $K_Y^2=9-\sum_{\nu} s_{\nu}$. Let $\bm{\delta}_{\nu}$ be the column vector of coefficients of $\pi^{\ast}K_S-K_Y$ along the exceptional components over the singular point indexed by $\nu$. By adjunction, we have $C_{\nu}\bm{\delta}_{\nu}=\bm{v}_{\nu}$. Thus, we have
\begin{equation}\label{eq:4.1}
    h=(-K_S)^2=9-\sum_{\nu}s_{\nu}+\sum_{\nu}\bm{v}_{\nu}^{\rm t}C_{\nu}^{-1}\bm{v}_{\nu}.
\end{equation}
For the chains $[3,2],[4,2]$, and $[3,2^5]$, the terms $\bm{v}_{\nu}^{\rm t}C_{\nu}^{-1}\bm{v}_{\nu}$ are $\frac{2}{5},\frac{8}{7}$, and $\frac{6}{13}$, respectively. For the chain $[2,3,2,2]$, we have $\bm{v}_{\nu}=(0,1,0,0)^{\rm t}$. Thus,
\[
\bm{v}_{\nu}^{\rm t}C_{\nu}^{-1}\bm{v}_{\nu}=(C_{\nu}^{-1})_{22}=\frac{6}{11}.
\] 
If the exceptional chain has type $[2^t]$, then $\bm v_{\nu}=0$. For a single smooth rational curve with self-intersection $-a$, we have $\bm{v}_{\nu}^{\rm t}C_{\nu}^{-1}\bm{v}_{\nu} =\frac{(a-2)^2}{a}$.

Let $D_{\rm max}\coloneqq \max_{\nu} \det C_{\nu}$. By \eqref{eq:4.1} and direct computation, we have
\begin{alignat*}{3}
\mathrm{I}\colon\quad &h=-1+2\cdot\frac25+\frac13=\frac2{15},
&\qquad&D_{\max}=5,\\
\mathrm{II}\colon\quad &h=-2+\frac87+\frac6{13}+\frac25=\frac2{455},
&&D_{\max}=13,\\
\mathrm{III}\colon\quad &h=-2+\frac87+\frac6{11}+\frac13=\frac5{231},
&&D_{\max}=11,\\
\mathrm{No.}\ 1\colon\quad &h=\frac13,
&&D_{\max}=5,\\
\mathrm{No.}\ 2\colon\quad &h=-1+\frac87=\frac17,
&&D_{\max}=7,\\
\mathrm{No.}\ 3\colon\quad &h=-2+\frac13+\frac95=\frac2{15},
&&D_{\max}=5.
\end{alignat*}
The values of $D_{\rm max}$ follow from $d([2^t])=t+1$, $d([3,2])=5$, $d([4,2])=d([2,4])=7$, $d([3,2^5])=13$, and $d([2,3,2,2])=11$. Note that every row satisfies $h<4$ and $hD_{\rm max}<4$. Thus, Corollary \ref{cor:2.6} excludes both a smooth base point and a singular base point for a pencil associated to a cylinder. Hence, no surface in rows I--III or Nos. 1--3 contains a cylinder.

\subsection{The birational construction}\label{subsec:4.3}
Let $T$ be the degree one Du Val del Pezzo surface with two points of type ${\sf A}_4$, and let $Y_0\to T$ be the minimal resolution. In characteristic five, $|-K_T|$ contains a cuspidal curve $\Gamma_T$ contained in the smooth locus of $T$. The strict transform $\Gamma\subset Y_0$ belongs to $|-K_{Y_0}|$ and satisfies $\Gamma^2=1$ (see \cite[Example 2.10]{Nag25}).

Blow up the cusp of $\Gamma$, then the tangency point of the strict transform of $\Gamma$ with the first exceptional curve, then the point where the strict transform of $\Gamma$ meets the two exceptional curves. Denote the resulting surface by $Y_3$ and the exceptional curve of the third blow-up by $F_0$. 

On $Y_3$, let $T_2$ be the strict transform of the exceptional curve of the second blow-up, let $T_3$ be the strict transform of the exceptional curve of the first blow-up, and let $T_5$ be the strict transform of $\Gamma$. Then $T_2^2=-2$, $T_3^2=-3$, $T_5^2=-5$, and $F_0^2=-1$. The curves $T_2,T_3,T_5$ are pairwise disjoint, and each meets $F_0$ transversely at a single point. These three intersection points are distinct. Moreover,
\begin{equation}\label{eq:4.2}
    -K_{Y_3}\sim T_2+T_3+T_5+2F_0.
\end{equation}

By \cite[Definition 5.5 and Lemmas 6.1--6.3]{Nag25}, the minimal resolution of every surface in Nos. 4--21 is obtained by further blow-ups on the divisor $T_2+T_3+T_5+F_0$, as described in the rest of this subsection and in Table \ref{tab:2}. At each stage, $Y$ denotes the resulting surface, $L_Y$ the exceptional curve of the most recent blow-up, and $F_0,T_2,T_3,T_5$ also denote their strict transforms on $Y$.

For No. 4, blow up a point of $F_0\setminus (T_2\cup T_3\cup T_5)$, then blow up $n$ times at the intersection of the strict transform of $F_0$ with $L_Y$. The resulting exceptional divisors are as follows:
\[
\begin{aligned}
 E^{(1)}\text{ has type }[n+2;[2],[3],[5]],\quad E^{(2)}\text{ has type }[2^n].
\end{aligned}
\]
Denote the central component of $E^{(1)}$, the strict transform of $F_0$, by $E_0$.  For $n=0$, we have $E^{(2)}=0$. The initial point can be chosen arbitrarily in $F_0\setminus (T_2\cup T_3\cup T_5)$, and the intersection numbers used below are independent of this choice (see also \cite[Theorem 1.2(2)]{Nag25}).

For Nos. 5--21, choose $r\in \{2,3,5\}$, and let $p<q$ be the other two integers in $\{2,3,5\}$. Blow up $F_0\cap T_r$, then blow up $u\geq 0$ times at the intersection of the strict transform of $F_0$ with the most recent exceptional curve. Denote the resulting surface by $Y_r(u)$. The two connected divisors joined by $L_Y$ are as follows:
\begin{equation}\label{eq:4.3}
\begin{aligned}
 E^{(1)}\text{ has type }[p,u+2,q],\quad E^{(2)}\text{ has type }[2^u,r+1].
\end{aligned}
\end{equation}
The divisor $E^{(1)}$ consists of the strict transforms of $T_p,F_0,T_q$, in that order. The chain $E^{(2)}$ is ordered starting from the component meeting $L_Y$, and the last component is the strict transform of $T_r$. 
The surfaces $Y_2(n)$, $Y_3(n)$, and $Y_5(n)$ are the minimal resolutions for Nos. 5, 6, and 7, respectively.

Starting from $Y_r(u)$, we obtain the minimal resolutions of the surfaces in Nos. 8--21 by performing the additional blow-ups listed in Table \ref{tab:2}. The index $i\in \{1,2\}$ specifies the center $L_Y\cap E^{(i)}$ of the blow-up. The new exceptional curve is denoted by $L_Y$, and the strict transform of the previous $L_Y$ is a $(-2)$-curve in $E^{(3-i)}$.

For Nos. 5--21, let $E_0$ be the strict transform of $F_0$, and write $E^{(1)}=T_p+E_0+T_q+B$, where $B=0$ for Nos. 5--7. We have $E_0^2=-(u+2)$. When $B\neq 0$, the chain $B$ joins $E_0$ to $L_Y$ and is ordered from $E_0$ toward $L_Y$. For Nos. 4--21, let $A=E^{(2)}$, ordered from the component meeting $L_Y$ whenever $A\neq 0$.

\begin{table}[htbp]
\centering
\small
\setlength{\tabcolsep}{6pt}
\renewcommand{\arraystretch}{1.25}
\begin{tabular}{@{}cccl@{}}
\toprule
Nos. & $r$ & Type of $B$ & Blow-ups after $Y_r(u)$\\
\midrule
$8,9$ & $2$ & $[2]$ & $2$\\
$10,11$ & $3$ & $[2^\ell]$ & $2$ repeated $\ell$ times\\
$12,13$ & $3$ & $[3]$ & $2,1$\\
$14,15$ & $5$ & $[2^\ell]$ & $2$ repeated $\ell$ times\\
$16,17$ & $5$ & $[\ell+2]$ & $2$, followed by $\ell$ indices equal to $1$\\
$18,19$ & $5$ & $[3,2]$ & $2,1,2$\\
$20,21$ & $5$ & $[2,3]$ & $2,2,1$\\
\bottomrule
\end{tabular}
\caption{Sequences of blow-ups for Nos. 8--21. In each row, $u=0$ corresponds to the even-numbered case, whereas $u=m\geq1$ corresponds to the odd-numbered case. The range of $\ell$ is specified in Table \ref{tab:1}.}
\label{tab:2}
\end{table}

For Nos. 5--21, write $A=A_1+\cdots+A_s$ and $B=B_1+\cdots+B_t$, with types $[a_1,\dots,a_s]$ and $[b_1,\dots,b_t]$, respectively. When $B=0$, the sequence $[b_1,\dots,b_t]$ is omitted. For $i=2$, blowing up at $L_Y\cap A_1$ changes the types to
\begin{equation}\label{eq:4.4}
\bigl([b_1,\ldots,b_t],[a_1,\ldots,a_s]\bigr)\longmapsto
\bigl([b_1,\ldots,b_t,2],\ [a_1+1,a_2,\ldots,a_s]\bigr).
\end{equation}
For $i=1$ and $B\neq 0$, blowing up $B_t\cap L_Y$ changes the types to
\begin{equation}\label{eq:4.5}
\bigl([b_1,\ldots,b_t],[a_1,\ldots,a_s]\bigr)\longmapsto
\bigl([b_1,\ldots,b_{t-1},b_t+1],\ [2,a_1,\ldots,a_s]\bigr).
\end{equation}

Let $\tau\colon Y\to Y_0$ be the resulting composition of blow-ups and $\pi\colon Y\to S$ the minimal resolution of a surface in Nos. 4--21. Then
\[
(\tau^*\Gamma)_{\mathrm{red}}=L_Y+E^{(1)}+E^{(2)}.
\]
The curve $L_Y$ is not contracted by $\pi$, whereas every component of $E^{(1)}+E^{(2)}$ is contracted by $\pi$, and the two exceptional chains of type ${\sf A}_4$ are disjoint from $(\tau^*\Gamma)_{\mathrm{red}}$. Define $P_1\coloneqq \pi(E^{(1)})$ and, when $E^{(2)}\neq 0$, $P_2\coloneqq \pi(E^{(2)})$. Then $\pi^{-1}(P_i)_{\rm red}=E^{(i)}$. The curve $L_Y$ meets $A_1$ when $A\neq 0$, and meets $E_0$ for Nos. 4--7 or $B_t$ for Nos. 8--21. All intersections are transverse of multiplicity one, and $L_Y$ meets no other exceptional component.

\subsection{Notation for the construction}
\begin{notation}\label{subsec:4.4}
     The surface $Y$, the curve $L_Y$, the exceptional divisors $E^{(1)}$, $E^{(2)}$, the points $P_1,P_2$, and the ordered chains $A,B$ are those constructed in Subsection \ref{subsec:4.3}. The point $P_2$ is used only when $E^{(2)}\neq0$, and the chain $B$ is used only for Nos. 5--21. Let $\pi, H, h, L, \eta, \delta_T$ and $\alpha_T$ be as in Subsection \ref{subsec:2.1}, defined for $S$ and the curve $L_Y$. Write $A=A_1+\cdots+A_s$ and $B=B_1+\cdots+ B_t$ when the chains are nonempty. The components of $A$ are ordered from $L_Y$ and the components of $B$ are ordered from $E_0$ toward $L_Y$.

    For each nonempty $E^{(i)}$, let $\mathcal{E}_i$ be the set of its irreducible components and write $C_i=C_{E^{(i)}}$, $G_i=G_{E^{(i)}}$, and $D_i=d(E^{(i)})$. The rows and columns are indexed by those components. Thus, $(G_i)_{T,T'}$ is the entry in row $T$ and column $T'$. If $E^{(2)}=0$, then $D_2=1$, $\mathcal{E}_2=\emptyset$. The remaining exceptional curves form two chains of type ${\sf A}_4$ and have coefficient zero in $\pi^{\ast}K_S-K_Y$.

    Let $N\coloneqq \frac{h}{\eta}$. If $A\neq 0$ has type $[a_1,\dots,a_s]$, then the determinants $p_i,q_i$ always refer to this ordered chain, as in \eqref{eq:2.3}, and $P_A\coloneqq p_s=d([a_1,\dots,a_{s-1}])$. If $B\neq 0$ has type $[b_1,\dots,b_t]$, then let
    \[
    d_B\coloneqq d(B),\quad \theta_B\coloneqq d([b_2,\dots,b_t]),\quad \xi_j\coloneqq d([b_1,\dots,b_{j-1}]),\quad \zeta_j\coloneqq d([b_{j+1},\dots,b_t]).
    \]
    Thus, $\xi_1=\zeta_t=1$ and $\zeta_1=\theta_B$. For $B=0$, let $d_B=1$ and $\theta_B=0$.

    The component of $E^{(1)}$ meeting $L_Y$ is denoted by $E_R$. We have $E_R=E_0$ for Nos. 4--7 and $E_R=B_t$ for Nos. 8--21. For Nos. 5--21, let the integers $r,p,q,u$ be as in Subsection \ref{subsec:4.3} and let
    \[
    \varepsilon\coloneqq pq-p-q,\quad \Delta\coloneqq u+2-\frac{1}{p}-\frac{1}{q}-\frac{\theta_B}{d_B},\quad \kappa\coloneqq \frac{1}{p}+\frac{1}{q}+\frac{1}{d_B}-1.
    \]
    For the estimates at $P_1$, we define
    \begin{equation}\label{eq:4.6}
    (c_0,\Lambda_1)\coloneqq
    \begin{cases}
    \left(6,N+\frac{7}{15}\right)&\text{for No.\ 4},\\[3pt]
    \left(p,\max\left\{d_B,P_A+\left(1-\frac{1}{q}\right)d_B\right\}\right)&\text{for Nos.\ 5--21}.
    \end{cases}
\end{equation}
\end{notation}

\subsection{The anticanonical divisor and the curve $L_Y$}
\begin{lemma}\label{lem:4.1}
    For a surface $S$ in Nos. 4--21, the number $N$ defined in Subsection \ref{subsec:4.4} is a positive integer such that $-K_S\sim NL$. If $A=E^{(2)}$ is nonempty, then $N=P_A+1$. For No. 4, $N=n+1$ for every $n\geq 0$.
\end{lemma}
\begin{proof}
    Set $\Xi_{Y_3}\coloneqq T_2+T_3+T_5+2F_0$, and for each subsequent blow-up $\sigma\colon Y\to V$ with exceptional curve $F$, define recursively $\Xi_Y\coloneqq \sigma^*\Xi_V - F$. Since $K_Y = \sigma^*K_V + F$, \eqref{eq:4.2} gives $\Xi_Y \sim -K_Y$ inductively. The divisor $\Xi_Y$ is supported on $L_Y+E^{(1)}+E^{(2)}$. For No. 4, the first blow-up is centered on $F_0$, whose coefficient is two, and every blow-up for Nos. 5--21 and every subsequent blow-up in No. 4 is centered at the intersection of two components with positive integer coefficients. Hence, the coefficient $N_0$ of $L_Y$ in $\Xi_Y$ is a positive integer. Since $\pi_*\Xi_Y=N_0L$, we have $-K_S\sim N_0L$. Intersecting with $-K_S$ gives $h=N_0\eta$. Hence, $N_0=N$ and $-K_S\sim NL$.

    Suppose that $A\neq 0$, and let $w_i$ be the coefficient of $A_i$ in $\Xi_Y$. Then $w_s=1$: for Nos. 5--21, the curve $A_s$ is the strict transform of $T_r$, whose coefficient is one, and for No. 4 with $n\geq 1$, the curve $A_s$ is the strict transform of the exceptional curve of the first blow-up, whose coefficient is also one. 

    Let $w_0=N$ and $w_{s+1}=0$. By adjunction and $\Xi_Y\sim -K_Y$, we have $a_iw_i-w_{i-1}-w_{i+1}=a_i-2$. For $y_i=w_i-1$, the equations become
    \[
    y_{i-1}=a_iy_i-y_{i+1},\quad y_s=0,\quad y_{s+1}=-1.
    \]
    By \eqref{eq:2.4}, we obtain $y_0=d([a_1,\dots,a_{s-1}])=P_A$. Since $y_0=N-1$, we have $N=P_A+1$.

    For No. 4 with $n\geq 1$, the chain $A$ of type $[2^n]$ gives $P_A=n$ and $N=n+1$. For $n=0$, the exceptional curve $L_Y$ of the blow-up at a smooth point of $F_0$ has coefficient one, since the coefficient of $F_0$ is two. Hence, $N=n+1$.
\end{proof}

\begin{lemma}\label{lem:4.2}
    For every surface $S$ in Nos. 4--21,
    \[
    L^2=\frac{1}{D_1D_2},\quad \eta=\frac{N}{D_1D_2},\quad h=\frac{N^2}{D_1D_2}.
    \]
\end{lemma}
\begin{proof}
    Let $(\tau^*\Gamma)_{\mathrm{red}}=\Gamma_Y+F_1+\cdots+F_m$, where $\Gamma_Y$ is the strict transform of $\Gamma$ and $F_i$ is the strict transform on $Y$ of the exceptional curve of the $i$th blow-up in $\tau$. Let $C_{\mathrm{all}}$ be the negative of the intersection matrix of $\Gamma_Y,F_1,\ldots,F_m$.

    For each blow-up, let $\widehat F_i$ denote the total transform on $Y$ of the exceptional curve created by the $i$th blow-up. The classes $\Gamma_Y,F_1,\ldots,F_m$ and those of $\tau^*\Gamma,\widehat F_1,\ldots,\widehat F_m$ form two bases of the same lattice. The change of basis matrix is integral and triangular with diagonal entries equal to $1$. Hence, the two intersection matrices have the same determinant.

    In the second basis, we have $(\tau^*\Gamma)^2=\Gamma^2=1$, $\tau^*\Gamma\cdot\widehat F_i=0$, and $\widehat F_i\cdot\widehat F_j=-\delta_{ij}$, where $\delta_{ij}$ is the Kronecker delta. Therefore, the intersection matrix in the basis $\tau^*\Gamma,\widehat F_1,\ldots,\widehat F_m$ is $\mathrm{diag}(1,-1,\ldots,-1)$. The negative of this diagonal matrix has determinant $(-1)\cdot 1^m=-1$. Consequently, $\det C_{\mathrm{all}}=-1$.

    For each nonempty $E^{(i)}=\sum_j E^{(i)}_j$, put $\mathbf d^{(i)}=(L_Y\cdot E^{(i)}_j)_j^{\rm t}$. Each vector $\mathbf d^{(i)}$ has one entry equal to one and every other entry zero. In the order $L_Y,E^{(1)},E^{(2)}$, the matrix $C_{\mathrm{all}}$ is
\[
C_{\mathrm{all}}=
\begin{pmatrix}
1&-(\mathbf d^{(1)})^{\rm t}&-(\mathbf d^{(2)})^{\rm t}\\
-\mathbf d^{(1)}&C_1&0\\
-\mathbf d^{(2)}&0&C_2
\end{pmatrix}.
\]
If $E^{(2)}=0$, the last block row and column are omitted. The Schur complement gives
\[
-1=D_1D_2\left(1-\sum_i(\mathbf d^{(i)})^{\rm t} C_i^{-1}\mathbf d^{(i)}\right),
\]
where the sum runs over the nonempty divisors $E^{(i)}$. Since $\pi^*L\cdot E_j^{(i)}=0$ for every exceptional component, the coefficient vector of $\pi^*L-L_Y$ along $E^{(i)}$ is
$C_i^{-1}\mathbf d^{(i)}$. Hence,
\[
L^2=-1+\sum_i(\mathbf d^{(i)})^{\rm t}C_i^{-1}\mathbf d^{(i)}.
\]
Thus, we obtain $L^2=\frac{1}{D_1D_2}$. Finally, $-K_S\sim NL$ implies $\eta=(-K_S)\cdot L=NL^2$ and $h=(-K_S)^2=N^2L^2$.
\end{proof}

\subsection{Determinants of the exceptional divisors}
For Nos. 5--21, we have
\[
(r,p,q,\varepsilon)=(2,3,5,7),\quad (3,2,5,3),\quad (5,2,3,1).
\]
In each case, we have
\begin{equation}\label{eq:4.7}
    r\varepsilon=pq-1,\quad p^2-\varepsilon\geq 1.
\end{equation}

\begin{lemma}\label{lem:4.3}
    For Nos. 5--21, the integers $P_A,N,D_1,D_2$ satisfy
    \begin{equation}\label{eq:4.8}
    \begin{aligned}
    P_A&=d_B(u+1)-\theta_B,&N&=P_A+1,\\
    D_2&=rP_A+d_B,&D_1&=\varepsilon D_2+P_A.
    \end{aligned}
    \end{equation}
    For No. 4, $N=D_2=n+1$ and $D_1=30n+29$.
\end{lemma}
\begin{proof}
    For an integer $b$, set $\mathsf M(b)\coloneqq\left(\begin{smallmatrix}b&-1\\1&0\end{smallmatrix}\right)$. For a nonzero ordered chain $A$ of type $[a_1,\ldots,a_s]$, define $\mathsf T(A)\coloneqq\mathsf M(a_1)\cdots\mathsf M(a_s)$. If $A=0$, define $\mathsf T(A)\coloneqq \mathrm{I}_2$. For $s\geq1$, equation \eqref{eq:2.4} gives
\[
\begin{aligned}
(\mathsf T(A))_{11}=d([a_1,\ldots,a_s]),\quad
(\mathsf T(A))_{12}=-d([a_1,\ldots,a_{s-1}]),\quad
(\mathsf T(A))_{21}=d([a_2,\ldots,a_s]).
\end{aligned}
\]
Each matrix $\mathsf T(A)$ has determinant one. Direct computation also gives
\[
\mathsf M(2)^u=\begin{pmatrix}u+1&-u\\u&1-u\end{pmatrix}\qquad(u\geq0).
\]

We first consider Nos. 5--21. The ordered chain $B+L_Y+A$ has matrix product $\mathsf T(B)\mathsf M(1)\mathsf T(A)$. The blow-ups in \eqref{eq:4.4} and \eqref{eq:4.5} preserve $\mathsf T(B)\mathsf M(1)\mathsf T(A)$. Indeed,
\begin{align*}
\mathsf M(2)\mathsf M(1)\mathsf M(b+1)=\mathsf M(1)\mathsf M(b),\quad
\mathsf M(b+1)\mathsf M(1)\mathsf M(2)=\mathsf M(b)\mathsf M(1).
\end{align*}
On $Y_r(u)$, one has $B=0$ and $A$ of type $[2^u,r+1]$. Consequently, at every subsequent stage of the construction, we have
\begin{equation}\label{eq:4.9}
\mathsf T(B)\mathsf M(1)\mathsf T(A)
=\begin{pmatrix}
r&-1\\r(u+1)+1&-(u+1)
\end{pmatrix}.
\end{equation}

The first column of $\mathsf T(B)$ is $(d_B,\theta_B)^{\rm t}$. Thus, the second row of $\mathsf T(B)^{-1}$ is $(-\theta_B,d_B)$. Since the first row of $\mathsf M(1)^{-1}$ is $(0,1)$, solving \eqref{eq:4.9} for $\mathsf T(A)$ gives
\begin{align*}
(D_2,-P_A)
&=(-\theta_B,d_B)\begin{pmatrix}r&-1\\r(u+1)+1&-(u+1)\end{pmatrix}\\
&=\bigl(r(d_B(u+1)-\theta_B)+d_B,\;-(d_B(u+1)-\theta_B)\bigr).
\end{align*}
Therefore, $P_A=d_B(u+1)-\theta_B$ and $D_2=rP_A+d_B$. The identity $N=P_A+1$ follows from Lemma \ref{lem:4.1}.

For nonempty $B$, Lemma \ref{lem:2.4} gives $(C_B^{-1})_{11}=\frac{\theta_B}{d_B}$. Taking the Schur complement of $\mathrm{diag}(p,q,C_B)$ in $C_1$, we obtain
\begin{equation}\label{eq:4.10}
\Delta= u+2-\frac1p-\frac1q-\frac{\theta_B}{d_B},
\qquad D_1=pqd_B\Delta.
\end{equation}
For $B=0$, eliminating $T_p,T_q$ gives 
\begin{equation}\label{eq:4.11}
    D_1=pq\left(u+2-\frac{1}{p}-\frac{1}{q}\right),
\end{equation}
which is \eqref{eq:4.10} with $d_B=1$ and $\theta_B=0$. Substituting $P_A=d_B(u+1)-\theta_B$ into $pqd_B\Delta$ yields
\[
D_1=pqP_A+(pq-p-q)d_B=pqP_A+\varepsilon d_B.
\]
Using $r\varepsilon=pq-1$ and $D_2=rP_A+d_B$, we obtain $D_1=\varepsilon D_2+P_A$.

For No. 4, $d([2^n])=n+1$ gives $D_2=n+1$. Taking the Schur complement of $\mathrm{diag}(2,3,5)$ in $C_1$ gives
$D_1=30(n+2-\frac{1}{2}-\frac{1}{3}-\frac{1}{5})=30n+29$.
Finally, Lemma \ref{lem:4.1} gives $N=n+1$, which completes the proof.
\end{proof}

\begin{lemma}\label{lem:4.4}
    For every surface in Nos. 4--21,
    \begin{equation}\label{eq:4.12}
        N\leq D_2<D_1,\quad 3N<D_1.
    \end{equation}
    For Nos. 5--21, $D_1<p^2D_2$ also holds.
\end{lemma}
\begin{proof}
    For Nos. 5--21, the integers $P_A,d_B,r,\varepsilon$ satisfy $P_A\geq1$, $d_B\geq1$, $r\geq2$, and $\varepsilon\geq1$. By \eqref{eq:4.8},
\[
D_2-N=(r-1)P_A+d_B-1\geq0,\qquad
D_1-D_2=(\varepsilon-1)D_2+P_A>0.
\]
By \eqref{eq:4.7} and \eqref{eq:4.8}, we have $D_1=pqP_A+\varepsilon d_B$. Since $pq\geq6$ and $N=P_A+1$, it follows that $D_1-3N=(pq-3)P_A+\varepsilon d_B-3\geq3P_A+d_B-3\geq1$. The inequality $p^2-\varepsilon\geq1$ in \eqref{eq:4.7} gives $p^2D_2-D_1=(p^2-\varepsilon)D_2-P_A\geq D_2-P_A>0$. For No. 4, the identities $N=n+1$, $D_2=N$, and $D_1=30N-1$ imply both inequalities in \eqref{eq:4.12} for $N\geq1$.
\end{proof}

\subsection{Excluding $P_2$}
\begin{proposition}\label{prop:4.5}
    Let $S$ be a surface in Nos. 4--21 with $E^{(2)}\neq 0$. Then $P_2=\pi(E^{(2)})$ is not the base point of any pencil associated to a cylinder in $S$.
\end{proposition}
\begin{proof}
    For the chain $A$ as in Subsection \ref{subsec:4.4}, Lemma \ref{lem:4.1} gives $p_s=P_A=N-1$, and $q_s=1$. The curve $L_Y$ meets $A_1$. Formula \eqref{eq:2.5} gives $(G_2)_{A_1,A_i}=\frac{q_i}{D_2}$ and $\frac{(G_2)_{A_i,A_i}}{(G_2)_{A_1,A_i}}=p_i$. Since $p_1<\cdots<p_s$,
\begin{equation}\label{eq:4.13}
\max_{1\leq i\leq s}\frac{(G_2)_{A_i,A_i}}{(G_2)_{A_1,A_i}}=P_A=N-1.
\end{equation}

For $s\geq2$, adjunction on $A_s$ and the equality $q_{s-1}=a_s q_s$ give $\delta_{A_{s-1}}q_s-\delta_{A_s}q_{s-1}=-(a_s-2)q_s$. For $2\leq j\leq s-1$, using $a_j\delta_{A_j}-\delta_{A_{j-1}}-\delta_{A_{j+1}}=a_j-2$ and $q_{j-1}=a_jq_j-q_{j+1}$, we obtain $\delta_{A_{j-1}}q_j-\delta_{A_j}q_{j-1}=\delta_{A_j}q_{j+1}-\delta_{A_{j+1}}q_j-(a_j-2)q_j$. It follows that 
\[
\delta_{A_i}q_{i+1}-\delta_{A_{i+1}}q_i=-\sum_{j=i+1}^s(a_j-2)q_j\leq0\qquad (1\leq i<s).
\]
Since $q_i>0$ for every $i$, we have 
\[
\frac{\delta_{A_i}}{q_i} \leq \frac{\delta_{A_{i+1}}}{q_{i+1}} \qquad (1\leq i<s).
\]
In particular, $\frac{\delta_{A_i}}{q_i}\leq\delta_{A_s}$, since $q_s=1$.
When $s=1$, the conclusion is immediate since $q_s=1$. Hence, by Lemma \ref{lem:2.4},
\begin{equation}\label{eq:4.14}
\max_{1\leq i\leq s}\frac{\delta_{A_i}}{(G_2)_{A_1,A_i}}
=D_2\delta_{A_s}=D_2-(p_s+q_s)=D_2-N.
\end{equation}

Since $\frac{h}{\eta}=N$, equations \eqref{eq:4.13} and \eqref{eq:4.14} give the two numbers in Proposition \ref{prop:2.8} explicitly:
\begin{equation}\label{eq:4.15}
\gamma_{P_2}=\frac{D_2}{N},\qquad
\beta_{P_2}=\frac{N-1}{D_1},\qquad
\gamma_{P_2}^2h=\frac{D_2}{D_1}.
\end{equation}
Lemma \ref{lem:4.4} gives $\beta_{P_2}<\frac{1}{3}$ and $\gamma_{P_2}^2h<1$. Therefore,
$\beta_{P_2}<1$ and $\gamma_{P_2}^2h<1<4(1-\beta_{P_2})^2$.
Proposition \ref{prop:2.8} excludes $P_2$ as the base point of a pencil associated to a cylinder.
\end{proof}

\subsection{The matrix and discrepancies over $P_1$}
For Nos. 5--21, the construction gives $E_0^2=-(u+2)$, $T_p^2=-p$, and $T_q^2=-q$. We now determine the two maxima appearing in Proposition \ref{prop:2.8}, using the component $E_R$ and the matrices defined in Subsection \ref{subsec:4.4}.

\begin{lemma}\label{lem:4.6}
For every surface in Nos. 5--21,
\begin{equation*}
\max_{D\in\mathcal E_1}\frac{(G_1)_{D,D}}{(G_1)_{E_R,D}}=\Lambda_1,
\qquad 2\Lambda_1\leq pD_2.
\end{equation*}
\end{lemma}
\begin{proof}
    Suppose first that $B\neq 0$. Let $\xi_j,\zeta_j$ be as in Subsection \ref{subsec:4.4}, and let $\bm{e}_1=(1,0,\dots,0)^{\rm t}\in\Q^t$ be the first coordinate vector. In the order $E_0,T_p,T_q,B_1,\ldots,B_t$, the matrix $C_1$ is
\begin{equation*}
C_1=\begin{pmatrix}
u+2&-1&-1&-\bm{e}_1 ^{\rm t}\\
-1&p&0&0\\
-1&0&q&0\\
-\bm{e}_1&0&0&C_B
\end{pmatrix}.
\end{equation*}
Lemma \ref{lem:2.4} gives
$(C_B^{-1})_{ij}=\frac{\xi_{\min\{i,j\}}\zeta_{\max\{i,j\}}}{d_B}$ and $C_B^{-1}\bm{e}_1=\left(\frac{\zeta_j}{d_B}\right)_{1\leq j\leq t} ^{\rm t}$. The Schur complement of $\mathrm{diag}(p,q,C_B)$ is $\Delta$ by \eqref{eq:4.10}.

The entries of $G_1=C_1^{-1}$ in the row indexed by $E_R=B_t$ are
\begin{equation}\label{eq:4.16}
\begin{aligned}
(G_1)_{E_R,E_0}&=\frac1{d_B\Delta},&
(G_1)_{E_R,T_p}&=\frac1{d_Bp\Delta},&
(G_1)_{E_R,T_q}&=\frac1{d_Bq\Delta},\\
(G_1)_{E_R,B_j}&=\frac{\xi_j}{d_B}+\frac{\zeta_j}{d_B^2\Delta}
&&\quad(1\leq j\leq t).
\end{aligned}
\end{equation}
The diagonal entries of $G_1$ are
\begin{equation}\label{eq:4.17}
\begin{aligned}
(G_1)_{E_0,E_0}&=\frac1\Delta,\\
(G_1)_{T_p,T_p}&=\frac1p+\frac1{p^2\Delta},
& (G_1)_{T_q,T_q}&=\frac1q+\frac1{q^2\Delta},\\
(G_1)_{B_j,B_j}&=\frac{\xi_j\zeta_j}{d_B}+\frac{\zeta_j^2}{d_B^2\Delta}
&&\quad(1\leq j\leq t).
\end{aligned}
\end{equation}
Dividing each diagonal entry in \eqref{eq:4.17} by the entry in row $E_R$ and the same column in \eqref{eq:4.16} gives
\begin{equation}\label{eq:4.18}
\begin{aligned}
\frac{(G_1)_{E_0,E_0}}{(G_1)_{E_R,E_0}}&=d_B,&
\frac{(G_1)_{B_j,B_j}}{(G_1)_{E_R,B_j}}&=\zeta_j\quad(1\leq j\leq t),\\
\frac{(G_1)_{T_p,T_p}}{(G_1)_{E_R,T_p}}&=d_B\left(\Delta+\frac1p\right),&
\frac{(G_1)_{T_q,T_q}}{(G_1)_{E_R,T_q}}&=d_B\left(\Delta+\frac1q\right).
\end{aligned}
\end{equation}
The inequalities $\zeta_j<d_B$ follow from Lemma \ref{lem:2.4}. Since $p<q$ and $d_B(\Delta+\frac{1}{p}) =P_A+(1-\frac{1}{q})d_B$, the maximum in \eqref{eq:4.18} is $\Lambda_1$.

If $B=0$, then $N=u+2$ by Lemma \ref{lem:4.1} and \eqref{eq:4.8}. In the order $T_p,E_0,T_q$, computing $C_1^{-1}$ gives
\[
G_1=\frac1{D_1}\begin{pmatrix}
qN-1&q&1\\q&pq&p\\1&p&pN-1
\end{pmatrix}.
\]
We have
$\frac{(G_1)_{T_p,T_p}}{(G_1)_{E_0,T_p}}=N-\frac{1}{q}$, $\frac{(G_1)_{E_0,E_0}}{(G_1)_{E_0,E_0}}=1$, and $\frac{(G_1)_{T_q,T_q}}{(G_1)_{E_0,T_q}}=N-\frac{1}{p}$. Since $P_A=N-1$, the maximum is again $\Lambda_1$.

We finally prove $2\Lambda_1\leq pD_2$. Since $D_2=rP_A+d_B$, we have
\begin{align*}
pD_2-2d_B&=prP_A+(p-2)d_B\geq0,\\ 
pD_2-2\left(P_A+\left(1-\frac1q\right)d_B\right) 
&=(pr-2)P_A+\left(p-2+\frac2q\right)d_B\geq0.
\end{align*}
Hence, we obtain $\Lambda_1\leq \frac{pD_2}{2}$.
\end{proof}

\begin{lemma}\label{lem:4.7}
    For every surface in Nos. 5--21,
    \begin{equation}\label{eq:4.19}
    \max_{D\in\mathcal E_1}\frac{\delta_D}{(G_1)_{E_R,D}}=\frac{D_1}{p}-N.
    \end{equation}
    The maximum is attained at $D=T_q$.
\end{lemma}
\begin{proof}
    Recall from Subsection \ref{subsec:4.4} that
\begin{equation}\label{eq:4.20}
\kappa=\frac1p+\frac1q+\frac1{d_B}-1.
\end{equation}
Suppose first that $B\neq 0$, and let $\bm{\alpha}_B=(\alpha_{B_1},\dots,\alpha_{B_t})^{\rm t}$. By adjunction, we obtain $p\alpha_{T_p}-\alpha_{E_0}=1$, $q\alpha_{T_q}-\alpha_{E_0}=1$, and $C_B\bm{\alpha}_B=\alpha_{E_0}\bm{e}_1+\bm{e}_t$, where $\bm{e}_1$ and $\bm{e}_t$ are the coordinate vectors of $\Q^t$. By Lemma \ref{lem:2.4}, $\alpha_{B_j}=\frac{\zeta_j\alpha_{E_0}+\xi_j}{d_B}$. The equation at $E_0$ is
\[
(u+2)\alpha_{E_0}-\alpha_{T_p}-\alpha_{T_q}-\alpha_{B_1}=-1.
\]
Thus, we obtain $\Delta \alpha_{E_0}=\kappa$. Hence,
\begin{equation}\label{eq:4.21}
\begin{aligned}
\alpha_{E_0}&=\frac\kappa\Delta,\\
\alpha_{T_p}&=\frac{1+\alpha_{E_0}}p,&
\alpha_{T_q}&=\frac{1+\alpha_{E_0}}q,\\
\alpha_{B_j}&=\frac{\zeta_j\alpha_{E_0}+\xi_j}{d_B}
&&\quad(1\leq j\leq t).
\end{aligned}
\end{equation}
If $B=0$, then the equation at $E_0$ is $(u+2)\alpha_{E_0}-\alpha_{T_p}-\alpha_{T_q}=0$. Together with the equations at $T_p$ and $T_q$, this gives the first three formulas in \eqref{eq:4.21} with $d_B=1$ and $\theta_B=0$.

Combining \eqref{eq:4.21} with the inverse matrix entries gives
\begin{equation*}
\begin{aligned}
\frac{\delta_{E_0}}{(G_1)_{E_R,E_0}}&=d_B(\Delta-\kappa),\\
\frac{\delta_{T_p}}{(G_1)_{E_R,T_p}}&=d_B((p-1)\Delta-\kappa),\\
\frac{\delta_{T_q}}{(G_1)_{E_R,T_q}}&=d_B((q-1)\Delta-\kappa).
\end{aligned}
\end{equation*}
Since $q>p\geq2$, the value at $T_q$ is at least the values at $E_0$ and $T_p$. If $B=0$, then $E^{(1)}=T_p+E_0+T_q$ and the maximum is attained at $T_q$.

Suppose that $B$ is nonempty, and write $\omega\coloneqq(q-1)\Delta-\kappa$. The coefficients $\delta_D$ are nonnegative since $\pi$ is the minimal resolution, and in particular, $\omega\geq0$. For each component $B_j$, formulas \eqref{eq:4.16} and \eqref{eq:4.21} give
\begin{equation*}
\frac{\delta_{B_j}}{(G_1)_{E_R,B_j}}
=d_B\,\frac{\Delta(d_B-\xi_j)-\zeta_j\kappa}
            {d_B\Delta\xi_j+\zeta_j}.
\end{equation*}
Since $d_B\Delta\xi_j+\zeta_j>0$, the inequality $\frac{\delta_{B_j}}{(G_1)_{E_R,B_j}}\leq d_B\omega$ is equivalent to
\begin{equation}\label{eq:4.22}
d_B\xi_j\omega+(q-1)\zeta_j-d_B+\xi_j\geq0.
\end{equation}
For $B$ of type $[2^\ell]$, the determinants are $\xi_j=j$, $\zeta_j=\ell+1-j$, and $d_B=\ell+1$. Thus, $d_B\xi_j\omega+(q-1)\zeta_j-d_B+\xi_j = d_B\xi_j\omega+(q-2)\zeta_j \geq0$.

For the chains $B$ in Table \ref{tab:2} which are not of type $[2^\ell]$, equations \eqref{eq:4.10} and \eqref{eq:4.20} give
\begin{alignat*}{2}
\text{type }[3],\ (p,q)=(2,5)&\colon\quad &\omega&=4u+\frac{23}{6},\\
\text{type }[b],\ (p,q)=(2,3),\ 3\leq b\leq5&\colon &\omega&=2u+\frac52-\frac3b,\\
\text{type }[3,2],\ (p,q)=(2,3)&\colon &\omega&=2u+\frac32,\\
\text{type }[2,3],\ (p,q)=(2,3)&\colon &\omega&=2u+\frac{11}{10}.
\end{alignat*}
For $B$ of type $[3]$, we have $(d_B,\theta_B)=(3,1)$; for $B$ of type $[b]$, $(d_B,\theta_B)=(b,1)$; for $B$ of type $[3,2]$, $(d_B,\theta_B)=(5,2)$; and for $B$ of type $[2,3]$, $(d_B,\theta_B)=(5,3)$. Each value of $\omega$ is at least one for $u\geq0$. Hence, in all four types,
\[
d_B\xi_j\omega+(q-1)\zeta_j-d_B+\xi_j
\geq d_B(\xi_j-1)+(q-1)\zeta_j+\xi_j>0.
\]
Inequality \eqref{eq:4.22} therefore holds for every component of $B$. Thus, the maximum in \eqref{eq:4.19} is attained at $T_q$.

Finally, $P_A=d_B(u+1)-\theta_B$ and $N=P_A+1$ give
$d_B(\Delta+\kappa)=d_B(u+1)-\theta_B+1=N$.
Since $D_1=pqd_B\Delta$, the value at $T_q$ becomes
\[
d_B((q-1)\Delta-\kappa)=qd_B\Delta-N=\frac{D_1}{p}-N.
\]
For $B=0$, the same formulas follow by setting $d_B=1$ and $\theta_B=0$.
\end{proof}

We next compute the two maxima in Proposition \ref{prop:2.8} for No. 4.

\begin{lemma}\label{lem:4.8}
    For No. 4, with the notation of Subsection \ref{subsec:4.4}, we have
\begin{equation*}
\max_{D\in\mathcal E_1}\frac{(G_1)_{D,D}}{(G_1)_{E_0,D}}=\Lambda_1,
\qquad
\max_{D\in\mathcal E_1}\frac{\delta_D}{(G_1)_{E_0,D}}=\frac{D_1}{6}-N.
\end{equation*}
Moreover, $2\Lambda_1\leq6D_2$ and $D_1<36D_2$.
\end{lemma}
\begin{proof}
    For No. 4, $D_1=30N-1$ and $D_2=N$. Order the components as $E_0,T_2,T_3,T_5$, and put $\boldsymbol\delta\coloneqq(\delta_{E_0},\delta_{T_2},\delta_{T_3},\delta_{T_5}) ^{\rm t}$. The inverse matrix $G_1$ and the vector of coefficients of $\pi^{\ast}K_S-K_Y$ are
\[
G_1=\frac1{30N-1}\begin{pmatrix}
30&15&10&6\\
15&15N+7&5&3\\
10&5&10N+3&2\\
6&3&2&6N+1
\end{pmatrix},\qquad
\boldsymbol\delta=\frac1{30N-1}\begin{pmatrix}
30N-2\\15N-1\\20N-1\\24N-1
\end{pmatrix}.
\]
Direct multiplication by $C_1$, whose diagonal entries are $N+1,2,3,5$ and whose off-diagonal entries are $-1$ in positions $(E_0,T_j)$ and $(T_j,E_0)$ for $j\in \{2,3,5\}$, gives $C_1G_1=\mathrm{I}_4$ and $C_1\boldsymbol\delta=(N-1,0,1,3) ^{\rm t}$.
Hence, we obtain
\begin{align*}
\frac{(G_1)_{E_0,E_0}}{(G_1)_{E_0,E_0}}&=1,&
\frac{(G_1)_{T_2,T_2}}{(G_1)_{E_0,T_2}}&=N+\frac7{15},\\
\frac{(G_1)_{T_3,T_3}}{(G_1)_{E_0,T_3}}&=N+\frac3{10},&
\frac{(G_1)_{T_5,T_5}}{(G_1)_{E_0,T_5}}&=N+\frac16.
\end{align*}
For the coefficients $\delta_D$, we have
\begin{align*}
\frac{\delta_{E_0}}{(G_1)_{E_0,E_0}}&=N-\frac1{15},&
\frac{\delta_{T_2}}{(G_1)_{E_0,T_2}}&=N-\frac1{15},\\
\frac{\delta_{T_3}}{(G_1)_{E_0,T_3}}&=2N-\frac1{10},&
\frac{\delta_{T_5}}{(G_1)_{E_0,T_5}}&=4N-\frac16.
\end{align*}
The first maximum is attained at $T_2$ and equals $\Lambda_1$. The second maximum is attained at $T_5$ and equals $4N-\frac{1}{6}=\frac{D_1}{6}-N$.
Finally, $6D_2-2\Lambda_1=4N-\frac{14}{15}>0$ and $36D_2-D_1=6N+1>0$ for $N\geq1$.
\end{proof}

\subsection{Excluding $P_1$}
\begin{proposition}\label{prop:4.9}
    For every surface $S$ in Nos. 4--21, the point $P_1=\pi(E^{(1)})$ is not the base point of any pencil associated to a cylinder in $S$.
\end{proposition}
\begin{proof}
    Let $E_R$ and $G_1$ be as in Subsection \ref{subsec:4.4}, and $c_0, \Lambda_1$ be as in \eqref{eq:4.6}. Lemmas \ref{lem:4.6}, \ref{lem:4.7}, and \ref{lem:4.8} give
\[
\max_{D\in\mathcal E_1}\frac{(G_1)_{D,D}}{(G_1)_{E_R,D}}=\Lambda_1,
\qquad
\max_{D\in\mathcal E_1}\frac{\delta_D}{(G_1)_{E_R,D}}=\frac{D_1}{c_0}-N,
\qquad 2\Lambda_1\leq c_0D_2.
\]
By Lemma \ref{lem:4.2}, the numbers defined in Proposition \ref{prop:2.8} satisfy
\begin{equation}\label{eq:4.23}
\gamma_{P_1}=\frac{D_1}{c_0N},\qquad
\beta_{P_1}=\frac{\Lambda_1}{c_0D_2}\leq\frac12,\qquad
\gamma_{P_1}^2h=\frac{D_1}{c_0^2D_2}.
\end{equation}
For Nos. 5--21, Lemma \ref{lem:4.4} gives $D_1<c_0^2D_2$ since $c_0=p$. For No. 4, Lemma \ref{lem:4.8} gives $D_1<c_0^2D_2$ since $c_0=6$. Thus,
$\beta_{P_1}<1$ and $\gamma_{P_1}^2h<1\leq4(1-\beta_{P_1})^2$.
Proposition \ref{prop:2.8} excludes $P_1$ as the base point of a pencil associated to a cylinder.
\end{proof}

\begin{proof}[Proof of the \ref{thm:main}]
    By \cite[Theorem 1.2 and Table 1]{Nag25}, every log del Pezzo surface of Picard rank one in characteristic five which is not log liftable appears in Table \ref{tab:1}. Subsection \ref{subsec:4.2} excludes cylinders on the surfaces in rows I--III and Nos. 1--3. For Nos. 4--21, the birational constructions are described in Subsection \ref{subsec:4.3}. 

    Let $S$ be one of the surfaces in Nos. 4--21, and suppose that $S$ contains a cylinder. By Lemma \ref{lem:2.1}, the associated pencil $\mathcal{L}$ has a base point $P$. Lemmas \ref{lem:4.2} and \ref{lem:4.4} give
    \[
    (-K_S)^2=\frac{N^2}{D_1D_2}\leq\frac N{D_1}<\frac13.
    \]
    By Corollary \ref{cor:2.6}, if $P$ is smooth, then $(-K_S)^2>4$, while if $P$ is of type $\mathsf A_4$, then $(-K_S)^2>\frac{4}{5}$. Thus, $P$ is either $P_1=\pi(E^{(1)})$ or, when $E^{(2)}$ is nonempty, $P_2=\pi(E^{(2)})$. Proposition \ref{prop:4.9} excludes $P_1$. When $E^{(2)}\neq0$, Proposition \ref{prop:4.5} also excludes $P_2$. Therefore, $S$ contains no cylinder.
\end{proof}

\section*{Acknowledgement}
This project was initiated during the 24th Affine Algebraic Geometry Meeting, held at Hirosaki University on March 4--6, 2026. The first two authors thank the organizers, Professors Adrien Dubouloz, Takashi Kishimoto, Hideo Kojima, and Masatomo Sawahara, for their warm hospitality. The first author is supported by a KIAS Individual Grant (MG111201) at Korea Institute for Advanced Study. The second author is supported by JSPS KAKENHI Grant Number JP21K13768. The third author is supported by JSPS KAKENHI Grant Numbers JP24K22823 and JP25K17222. 

\bibliographystyle{habbvr}
\bibliography{biblio}

\end{document}